\documentclass[a4, 12pt]{article}
\usepackage{authblk}
\usepackage{amsthm}
\usepackage{amsmath}
\usepackage{revsymb}

\makeatletter
\renewcommand*\env@matrix[1][*\c@MaxMatrixCols c]{%
  \hskip -\arraycolsep
  \let\@ifnextchar\new@ifnextchar
  \array{#1}}
\makeatother
\usepackage{mathptmx}
\usepackage{amssymb}
\usepackage{pdfpages}
\usepackage{graphics}
\usepackage{tikz}
\usepackage[margin=1in]{geometry}
\usepackage{setspace}

\newtheorem{theorem}{Theorem}

\newtheorem{corollary}{Corollary}
\newtheorem{case}{Case}

\numberwithin{subcase}{case}
\newtheorem{lemma}{Lemma}

\begin{document}
\title{Characterizing the Degree-Kirchhoff, Gutman, and Schultz Indices in Pentagonal Cylinders and M\"{o}bius Chains}

\author{Md. Abdus Sahir\thanks{The first author expresses gratitude to the University Grants Commission, India, for supporting this research through a Senior Research Fellowship. Email: \texttt{abdussahir@gmail.com}}~}

\author{Sk. Md. Abu Nayeem\thanks{Corresponding author. Email: \texttt{nayeem.math@aliah.ac.in}}}

\affil{Department of Mathematics and Statistics,\\ Aliah University, New Town, Kolkata -- 700160, India.}

\date{}
\maketitle

\begin{abstract} 
The degree-Kirchhoff index of a connected graph is defined as the sum of the reciprocals of the non-zero eigenvalues of the normalized Laplacian matrix, each multiplied by the graph's total degree. Several studies have recently obtained explicit formulations for the degree-Kirchhoff index of various kinds of class graphs. This paper presents closed-form formulas for the degree-Kirchhoff index of pentagonal cylinders and M\"{o}bius chains. Additionally, we calculate the Gutman index and Schultz index for these graphs.
\medskip



\end{abstract}

\maketitle

\section{Introduction}
Consider $G=(V,E)$, a connected, undirected graph with no self-loops or multiple edges. The vertices of $G$ are represented by $\{v_1,v_2,\ldots,v_{|V|}\}$, and its edges are represented by $E=\{e_1,e_2,\ldots,e_{|E|}\}$. The number of edges incident to a vertex $v_i$ in $G$ is its degree, denoted by $d_i$. The Laplacian matrix $L(G)=(l_{ij})$ and the normalized Laplacian matrix $\mathcal{L}(G)=(\ell_{ij})$ of the graph are defined as follows:
\[   
l_{ij}= 
     \begin{cases}
        -1  & \mbox{if } v_i\sim v_j\\
       d_i &\mbox{if } i=j\\
       0 &\mbox{otherwise,}\\
     \end{cases}
\]
\[   
\ell_{ij}= 
     \begin{cases}
        -\frac{1}{\sqrt{d_id_j}}  & \mbox{if } v_i\sim v_j\\
       1 &\mbox{if } i=j\\
       0 &\mbox{otherwise}\cdot\\
     \end{cases}
\]

The topological index is a numerical measure that characterizes the topological structure of a graph $G$. A topological index often includes a variety of graph parameters, including the order, size, distance between vertices pairs, degree of vertices, eccentricity, etc. The Wiener index \cite{Wie47}, a royal member of the family of old topological indices, was first presented by H. Wiener in $1947$ as $\operatorname{Wn}(G)=\frac{1}{2}\sum\limits_{i=1}^{|V|}\sum\limits_{j=1}^{|V|} d_{ij}$, where $d_{ij}$ represents the distance between the vertex pair $v_i$ and $v_j$. The Schultz index \cite{Sch89}, introduced by H. P. Schultz in $1989$ as $\operatorname{Sc}(G)=\frac{1}{2}\sum\limits_{i=1}^{|V|}\sum\limits_{j=1}^{|V|} (d_i+d_j)d_{ij},$ is a weighted version of the Wiener index. In $1994,$ another weighted variation of the Wiener index was presented by Gutman as  $\operatorname{Gut}(G)=\frac{1}{2}\sum\limits_{i=1}^{|V|}\sum\limits_{j=1}^{|V|} d_id_jd_{ij}$ which later came to be known as the Gutman index \cite{Gut94}.

Being motivated by the success of the Wiener index, nearly after four decades, in $1993,$ Klein and Randi\'{c} put forward a novel structure-descriptor topological index, known as the Kirchhoff index \cite{Kle93}. The Kirchhoff index of a graph $G$ is calculated using the formula $\operatorname{Kf}(G)=\frac{1}{2}\sum\limits_{i=1}^{|V|}\sum\limits_{j=1}^{|V|} r_{ij}$, where 
$r_{ij}$ represents the effective resistance between vertices 
$v_i$ and $v_j$. This effective resistance is determined by treating the edges as unit resistors and applying Ohm's law. In $2007,$ Chen and Zhang introduced the notion of the degree-Kirchhoff index \cite{Che07}, defined as $\operatorname{Kf}^*(G)=\frac{1}{2}\sum\limits_{i=1}^{|V|}\sum\limits_{j=1}^{|V|} d_id_jr_{ij}.$ For a connected graph, all eigenvalues of the Laplacian matrix and normalized Laplacian matrix are positive except the smallest one. Let  $\{0=\lambda_1< \lambda_2\leq \cdots\leq \lambda_{|V|}\},$ and $\{0=\lambdabar_1< \lambdabar_2\leq \cdots\leq \lambdabar_{|V|}\}$   denote the spectra of the Laplacian and normalized Laplacian matrices of $G$, respectively, where $(|V|\geq 2)$. Researchers have established a clear connection between the Kirchhoff index and the Laplacian spectrum, as well as between the degree-Kirchhoff index and the normalized Laplacian spectrum.

The Kirchhoff index can be interpreted alternatively as the following theorem, as demonstrated separately by Gutman and Mohar \cite{Gut96} and Zhu et al. \cite{Zhu96}.
\begin{theorem}{\cite{Gut96,Zhu96}}\label{Kiralt}
Let $G$ be a simple and connected graph. The Kirchhoff index of $G$ is 
\begin{eqnarray*}
\operatorname{Kf}(G)=|V|\sum\limits_{i=2}^{|V|}\frac{1}{\lambda_i}\cdot
\end{eqnarray*} 
\end{theorem}

Similarly, the degree-Kirchhoff index can be computed as follows.
\begin{theorem}{ \cite{Che07}}\label{degKiralt}
Let $G$ be a simple and connected graph. Then the degree-Kirchhoff index of $G$ is
\begin{eqnarray*}
\operatorname{Kf}^*(G)=2|E|\sum\limits_{i=2}^{|V|}\frac{1}{\lambdabar_i}\cdot
\end{eqnarray*} 
\end{theorem}

Kemeny's constant is a crucial graph invariant related to the study of the Markov chain  \cite{Lev02} that can also be calculated using the normalized Laplacian spectrum.

\begin{theorem} \cite{But16} \label{kemeny}
Kemeny's constant of a connected graph $G$ is 
\begin{eqnarray*}\operatorname{Kc}(G)=\sum\limits_{i=2}^{|V|}\frac{1}{\lambdabar_i},\end{eqnarray*} 
i.e., $\operatorname{Kc}(G)=\frac{1}{2|E|}\cdot \operatorname{Kf}^*(G)\cdot$ 
\end{theorem}

The following theorem is used to determine the total number of spanning trees in a connected graph $G$ using the Laplacian or normalized Laplacian spectrum \cite{Chu97}.  

\begin{theorem}{\cite{Chu97}}\label{spanningtree}
The total count of spanning trees of a connected graph $G$ is
\begin{eqnarray*} \tau(G)=\frac{1}{|V|} \prod\limits_{i=2}^{|V|} \lambda_i=\frac{1}{2|E|}\prod\limits_{i=1}^{|V|} d_i \prod\limits_{i=2}^{|V|} \lambdabar_i\cdot 
\end{eqnarray*}
\end{theorem}

As both of the structural descriptors have important applications in graph theory, networking systems, molecular chemistry, and other related fields, researchers have taken a strong interest in the Kirchhoff indices and degree-Kirchhoff indices of graphs in recent years. Finding the explicit formulae for the Kirchhoff index and the degree-Kirchhoff index for any general class of graphs is not an easy task. But, if a suitable automorphism is found for a particular class of graphs, then using the decomposition theorem, one can derive the explicit formula for them. 

Klein and Randi\'{c} \cite{Kle93} showed that the Kirchhoff index can not exceed the Wiener index and the equality holds only when $G$ is a tree. 
Gutman \cite{Gut94} demonstrated that for a tree $T$, the Wiener index and the Gutman index are related by the equation $\operatorname{Gut}(T)=4\operatorname{Wn}(T)-(2|V|-1)(|V|-1)$. Similarly, Klein et al. \cite{Kle92} established that the Wiener index and the Schultz index are connected by the formula $\operatorname{Sc}(T)=4\operatorname{Wn}(T)-|V|(|V|-1)$.

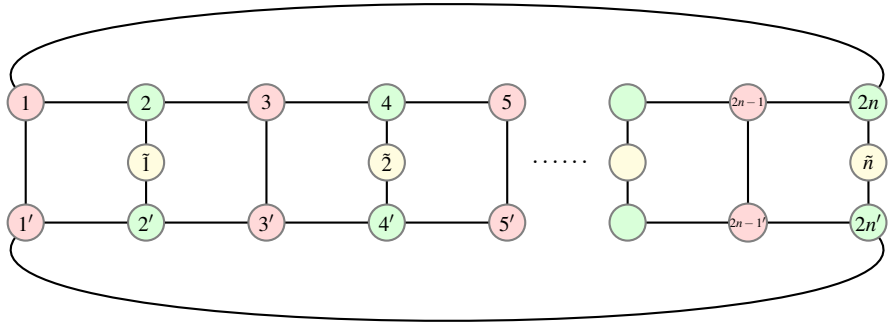
\begin{figure}
\begin{center}
\begin{tikzpicture}[place/.style={thick, circle,draw=black!50,fill=black!20,inner sep=0pt, minimum size = 6 mm}, transform shape, scale=0.8]
\node[place, fill=red!15!] (1) {\footnotesize 1} ;
\node (2)[place, fill=red!15! ] [below of = 1, node distance=2cm] {\footnotesize $1'$} edge[thick] (1);
\node (3)[place, fill=green!15!] [right of = 1, node distance=2cm] {\footnotesize 2} edge[thick] (1);

\node (4)[place, fill=yellow!15!] [below of = 3, node distance=1cm] {\footnotesize $\tilde{1}$} edge[thick] (3) ;

\node (5)[place, fill=green!15!] [ below of = 4] {\footnotesize $2'$} edge[thick] (4) edge[thick] (2);

\node (6)[place, fill=red!15!] [right of = 3,node distance=2cm] {\footnotesize 3} edge[thick] (3);
\node (7)[place, fill=red!15!] [below of = 6,node distance=2cm] {\footnotesize $3'$} edge[thick] (6) edge[thick] (5);
\node (8)[place, fill=green!15!] [right of = 6,node distance=2cm] {\footnotesize 4} edge[thick] (6);
\node (9)[place, fill=yellow!15!] [below of = 8, node distance=1cm] {\footnotesize $\tilde{2}$} edge[thick] (8) ;

\node (10)[place, fill=green!15!] [right of = 7,node distance=2cm] {\footnotesize $4'$} edge[thick] (7) edge[thick] (9);

\node (11)[place, fill=red!15!] [right of = 8,node distance=2cm] {\footnotesize 5} edge[thick] (8);

\node (12)[place, fill=red!15!] [below of = 11,node distance=2cm] {\footnotesize $5'$} edge[thick] (10) edge[thick] (11);

\node (13)[place, fill=green!15!] [right of = 11, node distance=2cm] {} ;
\node (14)[place, fill=yellow!15!] [below of = 13] {} edge[thick] (13);
\node (15)[place, fill=green!15!] [below of = 14, node distance=1cm] {}  edge[thick] (14);
\node (16)[place, fill=red!15!] [right of = 13, node distance=2cm] {\footnotesize \tiny $2n-1$} edge[thick] (13);
\node (17)[place, fill=red!15!] [below of = 16, node distance=2cm] {\footnotesize \tiny $2n-1'$} edge[thick] (16) edge[thick] (15);
\node (18)[place, fill=green!15!] [right of = 16, node distance=2cm] {\footnotesize $2n$} edge[thick] (16)  edge[thick, out=55, in=125,out looseness=0.4, in looseness=.4] (1);
\node (19)[place, fill=yellow!15!] [below of = 18, node distance=1cm] {\footnotesize $\tilde{n}$} edge[thick] (18);
\node (21)[place, fill=green!15!] [below of = 19, node distance=1cm] {\footnotesize $2n'$} edge[thick] (19) edge[thick] (17) edge[thick, out=-55, in=-125, out looseness=0.4, in looseness=0.4] (2);
\draw (11) node [xshift = 0.9cm, yshift =-1 cm] {$\cdots\cdots$};
\end{tikzpicture}
\caption{\label{fig1}  The pentagonal cylinder chain $P_n.$}
\end{center}
\end{figure}

\begin{figure}
\begin{center}
\begin{tikzpicture}[place/.style={thick, circle,draw=black!50,fill=black!20,inner sep=0pt, minimum size = 6 mm}, transform shape, scale=0.8]
\node [place, fill=red!15!] (1) {\footnotesize 1} ;
\node (2)[place, fill=red!15! ] [below of = 1, node distance=2cm] {\footnotesize $1'$} edge[thick] (1);
\node (3)[place, fill=green!15!] [right of = 1, node distance=2cm] {\footnotesize 2} edge[thick] (1);

\node (4)[place, fill=yellow!15!] [below of = 3, node distance=1cm] {\footnotesize $\tilde{1}$} edge[thick] (3) ;

\node (5)[place, fill=green!15!] [ below of = 4] {\footnotesize $2'$} edge[thick] (4) edge[thick] (2);

\node (6)[place, fill=red!15!] [right of = 3,node distance=2cm] {\footnotesize 3} edge[thick] (3);
\node (7)[place, fill=red!15!] [below of = 6,node distance=2cm] {\footnotesize $3'$} edge[thick] (6) edge[thick] (5);
\node (8)[place, fill=green!15!] [right of = 6,node distance=2cm] {\footnotesize 4} edge[thick] (6);
\node (9)[place, fill=yellow!15!] [below of = 8, node distance=1cm] {\footnotesize $\tilde{2}$} edge[thick] (8) ;

\node (10)[place, fill=green!15!] [right of = 7,node distance=2cm] {\footnotesize $4'$} edge[thick] (7) edge[thick] (9);

\node (11)[place, fill=red!15!] [right of = 8,node distance=2cm] {\footnotesize 5} edge[thick] (8);

\node (12)[place, fill=red!15!] [below of = 11,node distance=2cm] {\footnotesize $5'$} edge[thick] (10) edge[thick] (11);

\node (13)[place, fill=green!15!] [right of = 11, node distance=2cm] {} ;
\node (14)[place, fill=yellow!15!] [below of = 13] {} edge[thick] (13);
\node (15)[place, fill=green!15!] [below of = 14, node distance=1cm] {}  edge[thick] (14);
\node (16)[place, fill=red!15!] [right of = 13, node distance=2cm] {\footnotesize \tiny $2n-1$} edge[thick] (13);
\node (17)[place, fill=red!15!] [below of = 16, node distance=2cm] {\footnotesize \tiny $2n-1'$} edge[thick] (16) edge[thick] (15);
\node (18)[place, fill=green!15!] [right of = 16, node distance=2cm] {\footnotesize $2n$} edge[thick] (16)  edge[thick, out=55, in=125,out looseness=0.4, in looseness=1.2] (2);

\node (19)[place, fill=yellow!15!] [below of = 18, node distance=1cm] {\footnotesize $\tilde{n}$} edge[thick] (18);
\node (20)[place, fill=green!15!] [ below of = 19] {\footnotesize $2n'$} edge[thick] (19) edge[thick] (17) edge[thick, out=-55, in=-125, out looseness=0.4, in looseness=1.2] (1);
\draw (11) node [xshift = 0.9cm, yshift =-1 cm]{$\cdots\cdots$};

\end{tikzpicture}
\caption{\label{fig2}  The pentagonal M\"{o}bius chain $P'_n.$}
\end{center}
\end{figure}
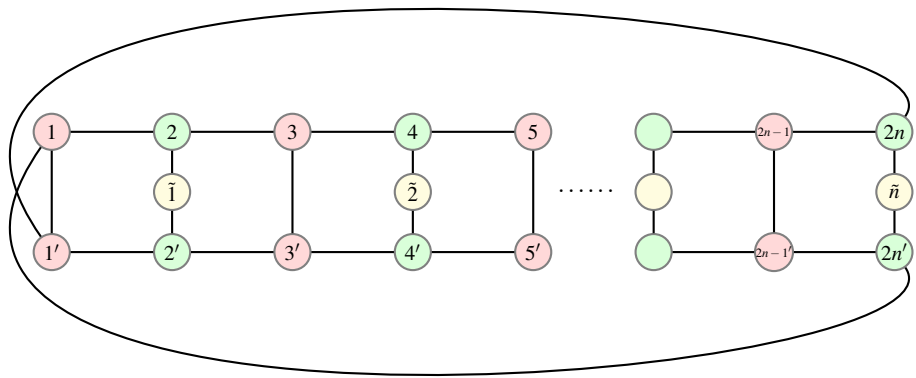

The formulae of the degree-Kirchhoff index of linear hexagonal chains \cite{Hua16}, hexagonal M\"{o}bius chains \cite{Ma19}, linear crossed hexagonal chains \cite{Pan18}, linear polyomino chains \cite{HuaJ16}, cylinder phenylene chains \cite{Ma21}, random polygonal chains \cite{Li21}, generalized phenylenes \cite{Zhu19}, cylinder, and  M\"{o}bius octagonal chain \cite{Liu22}, linear pentagonal derivation chain \cite{Liu24, Liu23} are already available in the literature. The formula to obtain the Kirchhoff index for pentagonal chains was obtained by Wang and Zhang \cite{Wan10} in $2010$ and the degree-Kirchhoff index for the same chain graphs was obtained by He et al. \cite{He18} in $2018.$ Although the formulae to find the Kirchhoff and degree-Kirchhoff indices of several cylinder/ M\"{o}bius chains have been computed by many researchers in the last few years, those for cylinder/ M\"{o}bius pentagonal chains have not been attempted. Very recently, explicit formulae for the Kirchhoff index of the cylinder pentagonal chains and the M\"{o}bius pentagonal chains have been obtained by Sahir and Nayeem \cite{Sah22}. We now present explicit formulae for the degree-Kirchhoff index, Kemeny's constant, Gutman index, and Schultz index of the pentagonal cylinder chain $P_n$ (see Figure \ref{fig1}) and the pentagonal M\"{o}bius chain $P'_n$ (see Figure \ref{fig2}) on $|V|=5n$ ($n \geq 2$) vertices and $|E|=7n$ edges, as a continuation of previous works in that direction. For the above-said graphs, we also derive a relationship between the Schultz and Gutman indices and present a comparison between the Gutman indices and the degree-Kirchhoff indices for different values of $n$.

\section{Preliminaries} An automorphism $\pi:V \rightarrow V$ is a bijective map with the properties that $u$ is adjacent to $v$ in $G$ if and only if $\pi(u)$ is adjacent to $\pi(v)$ in $G$. Let the automorphism $\pi$ gives the partition $V_0=\{v_{\tilde{1}},v_{\tilde{2}},\ldots,v_{\tilde{p}}\}$, $V_1=\{v_1,v_2,\ldots,v_q\}$ and $V_2=\{v_1',v_2',\ldots,v_q'\}$ of the vertex set $V$ where $\pi(v_{\tilde{i}})=v_{\tilde{i}}\mbox{ for all } v_{\tilde{i}}\in V_0, \pi(v_i)=v_i'\mbox{ for all }v_i\in V_1$ and $\pi(v_i')=v_i\mbox{ for all }v_i'\in V_2$. One can readily decompose $\pi$ as a product of disjoint permutation cycles of length $1$ and $2$, i.e., 
$$\pi=(v_{\tilde{1}})(v_{\tilde{2}})\cdots(v_{\tilde{p}})(v_1,v_1')(v_2,v_2')\cdots(v_q,v_q'),$$ where $p+2q=|V|$. The normalized Laplacian matrix $\mathcal{L}(G)$ can be expressed as follows once the vertices have been appropriately relabeled: 
$$\mathcal{L}(G)=\begin{pmatrix}
\mathcal{L}_{00} & \mathcal{L}_{01} & \mathcal{L}_{02}\\
\mathcal{L}_{10} & \mathcal{L}_{11} & \mathcal{L}_{12}\\
\mathcal{L}_{20} & \mathcal{L}_{21} & \mathcal{L}_{22}
\end{pmatrix},$$ where the submatrix $\mathcal{L}_{ij}$ corresponds to the vertices of the respective vertex sets $V_i$ and $V_j,i,j=0,1,2.$ 

Let, 
\begin{equation}\label{LA}\mathcal{L}_A(G)=\begin{pmatrix}[c|c]
\mathcal{L}_{00} & \sqrt{2}\mathcal{L}_{01} \\\hline
\sqrt{2}\mathcal{L}_{10} & \mathcal{L}_{11} + \mathcal{L}_{12}\\
\end{pmatrix}
\end{equation} and 
\begin{equation}\label{LS}\mathcal{L}_S(G)=\mathcal{L}_{11} - \mathcal{L}_{12}.\end{equation}

The following lemma is called the normalized Laplacian decomposition lemma.   
\begin{lemma}\label{Hua}\cite{Hua16}
Let $\phi(\mathcal{L}(G))=\det\left(yI-\mathcal{L}(G)\right)$ denotes the characteristic polynomial of the normalized Laplacian matrix $\mathcal{L}(G)$ where $I$ is the identity matrix of suitable size. Then $$\phi(\mathcal{L}(G))=\phi(\mathcal{L}_A(G))\cdot\phi(\mathcal{L}_S(G)),$$ where $\mathcal{L}_A(G)$ and $\mathcal{L}_S(G)$ are given by equations (\ref{LA}) and (\ref{LS}) respectively.
\end{lemma}

From now on, we use  $\mathcal{L}_A$ instead of $\mathcal{L}_A(P_n)$ and, $\mathcal{L}_S$ instead of $\mathcal{L}_S(P_n)$ for the sake of simplicity. Similar notations with primes ($'$) are used for pentagonal M\"{o}bius chain $P'_n$. Thus,

$\mathcal{L}_{00}=\mathcal{L}'_{00}=I_n$,
$\mathcal{L}_{01}=\mathcal{L}'_{01}=\begin{pmatrix}
0 & -\frac{1}{\sqrt{6}} &0&0& \cdots &0\\
0 & 0&0&-\frac{1}{\sqrt{6}}& \cdots &0\\
\vdots & \vdots & \vdots&\vdots  &\ddots &  \vdots\\
0 & 0&0&0& \cdots &-\frac{1}{\sqrt{6}}\\
\end{pmatrix}_{n\times 2n},$

$\mathcal{L}_{11}=\begin{pmatrix}
1 &-\frac{1}{3}&0& \cdots&0&-\frac{1}{3}\\
-\frac{1}{3} &1&-\frac{1}{3}& \cdots&0&0\\
0 & -\frac{1}{3}&1& \cdots &0&0\\
\vdots & \vdots & \vdots&\ddots   & \vdots& \vdots\\
0 &0&0& \cdots&1&-\frac{1}{3}\\
-\frac{1}{3} &0&0& \cdots&-\frac{1}{3}&1
\end{pmatrix}_{2n\times 2n}$,
$\mathcal{L}_{12}=\begin{pmatrix}
-\frac{1}{3} &0&0& \cdots&0&0\\
0 &0&0& \cdots&0&0\\
0 & 0&-\frac{1}{3}& \cdots &0&0\\
\vdots & \vdots & \vdots&\ddots   & \vdots& \vdots\\
0 &0&0& \cdots&-\frac{1}{3}&0\\
0 &0&0& \cdots&0&0
\end{pmatrix}_{2n\times 2n}$,

and $~\mathcal{L}'_{11}=\begin{pmatrix}
1 &-\frac{1}{3}&0& \cdots&0&0\\
-\frac{1}{3} &1&-\frac{1}{3}& \cdots&0&0\\
0 & -\frac{1}{3}&1& \cdots &0&0\\
\vdots & \vdots & \vdots&\ddots   & \vdots& \vdots\\
0 &0&0& \cdots&1&-\frac{1}{3}\\
0 &0&0& \cdots&-\frac{1}{3}&1
\end{pmatrix}_{2n\times 2n}$,
$\mathcal{L}'_{12}=\begin{pmatrix}
-\frac{1}{3} &0&0& \cdots&0&-\frac{1}{3}\\
0 &0&0& \cdots&0&0\\
0 & 0&-\frac{1}{3}& \cdots &0&0\\
\vdots & \vdots & \vdots&\ddots   & \vdots& \vdots\\
0 &0&0& \cdots&-\frac{1}{3}&0\\
-\frac{1}{3} &0&0& \cdots&0&0
\end{pmatrix}_{2n\times 2n}$.\\

So, 
$$\mathcal{L}_A=\mathcal{L}_A'=\begin{pmatrix}[cccc|cccccc]
1 & 0 & \cdots & 0  & 0 & -\frac{1}{\sqrt{3}} & 0 & 0 & \cdots & 0\\
0 &1&\cdots&0  &0 &0&0&-\frac{1}{\sqrt{3}}&\cdots&0 \\
\vdots&\vdots&\ddots&\vdots&\vdots&\vdots&\vdots&\vdots&\ddots&\vdots\\
0 & 0 & \cdots & 1  & 0 & 0 & 0 & 0 & \cdots & -\frac{1}{\sqrt{3}}\\\hline

0 & 0 & \cdots & 0           & \frac{2}{3}& -\frac{1}{3} & 0 & 0 & \cdots & -\frac{1}{3}\\
-\frac{1}{\sqrt{3}} & 0 & \cdots & 0    & -\frac{1}{3}& 1 & -\frac{1}{3} & 0 & \cdots & 0\\
0 & 0 & \cdots & 0              & 0& -\frac{1}{3} &\frac{2}{3} & -\frac{1}{3} & \cdots & 0\\
0 & -\frac{1}{\sqrt{3}}& \cdots & 0    & 0& 0& -\frac{1}{3} & 1 & \cdots & 0\\
\vdots&\vdots&\ddots&\vdots&\vdots&\vdots&\vdots&\vdots&\ddots&\vdots\\
0 & 0& \cdots & -\frac{1}{\sqrt{3}}   & -\frac{1}{3}& 0& 0 & 0 & \cdots & 1\\
\end{pmatrix}_{3n\times 3n},$$

\begin{eqnarray*}\mathcal{L}_S&=&\begin{pmatrix}
\frac{4}{3} &-\frac{1}{3}&0&0& \cdots&0&0&-\frac{1}{3}\\
-\frac{1}{3}&1 &-\frac{1}{3}&0& \cdots&0&0&0\\
0 &-\frac{1}{3}&\frac{4}{3}&-\frac{1}{3}& \cdots&0&0&0\\
0&0&-\frac{1}{3}&1 & \cdots&0&0&0\\
\vdots & \vdots & \vdots   & \vdots& \ddots& \vdots & \vdots & \vdots \\
0&0&0&0 & \cdots&1&-\frac{1}{3}&0\\
0&0&0&0 & \cdots&-\frac{1}{3}&\frac{4}{3}&-\frac{1}{3}\\
-\frac{1}{3}&0&0&0 & \cdots&0&-\frac{1}{3}&1\\
\end{pmatrix}_{2n\times 2n}
\end{eqnarray*}

\begin{eqnarray*}&=&\frac{1}{3}\cdot
\begin{pmatrix}
4 &-1&0&0& \cdots&0&0&-1\\
-1&3 &-1&0& \cdots&0&0&0\\
0 &-1&4&-1& \cdots&0&0&0\\
0&0&-1&3 & \cdots&0&0&0\\
\vdots & \vdots & \vdots   & \vdots& \ddots& \vdots & \vdots & \vdots \\
0&0&0&0 & \cdots&3&-1&0\\
0&0&0&0 & \cdots&-1&4&-1\\
-1&0&0&0 & \cdots&0&-1&3\\
\end{pmatrix}_{2n\times 2n},\end{eqnarray*}
and
\begin{eqnarray*}\mathcal{L}_S'&=&\begin{pmatrix}
\frac{4}{3} &-\frac{1}{3}&0&0& \cdots&0&0&\frac{1}{3}\\
-\frac{1}{3}&1 &-\frac{1}{3}&0& \cdots&0&0&0\\
0 &-\frac{1}{3}&\frac{4}{3}&-\frac{1}{3}& \cdots&0&0&0\\
0&0&-\frac{1}{3}&1 & \cdots&0&0&0\\
\vdots & \vdots & \vdots   & \vdots& \ddots& \vdots & \vdots & \vdots \\
0&0&0&0 & \cdots&1&-\frac{1}{3}&0\\
0&0&0&0 & \cdots&-\frac{1}{3}&\frac{4}{3}&-\frac{1}{3}\\
\frac{1}{3}&0&0&0 & \cdots&0&-\frac{1}{3}&1\\
\end{pmatrix}_{2n\times 2n}\\&=&\frac{1}{3}\cdot\begin{pmatrix}
4 &-1&0&0& \cdots&0&0&1\\
-1&3 &-1&0& \cdots&0&0&0\\
0 &-1&4&-1& \cdots&0&0&0\\
0&0&-1&3 & \cdots&0&0&0\\
\vdots & \vdots &\vdots   &\vdots & \ddots&  \vdots & \vdots & \vdots \\
0&0&0&0 & \cdots&3&-1&0\\
0&0&0&0 & \cdots&-1&4&-1\\
1&0&0&0 & \cdots&0&-1&3
\end{pmatrix}_{2n\times 2n}.\end{eqnarray*}

Let 
\begin{eqnarray}&&\phi(\mathcal{L}_A)=y^{3n}+\gamma_1y^{3n-1}+\cdots+\gamma_{3n-2}y^{2}+\gamma_{3n-1}y ~~(\mbox{since } \mathcal{L}_A \mbox{ is singular},\gamma_{3n}=0)\label{char1},\\ &&\mbox{and }\phi(\mathcal{L}_S)=y^{2n}+\delta_1y^{2n-1}+\cdots+\delta_{2n-2}y^{2}+\delta_{2n-1}y+\delta_{2n}\label{char2}.\end{eqnarray} 

Therefore, $\gamma_i=(-1)^i\cdot$ ($\sum\limits_{i=1}^{3n}$principal minors of $\mathcal{L}_A$ of order $i)$, and $\delta_j=(-1)^j\cdot$ ($\sum\limits_{j=1}^{2n}$principal minors of $\mathcal{L}_S$ of order $j)$.

Let $\{0=\rho_1<\rho_2\le\cdots\le\rho_{3n}\}$ represent the eigenvalues of $\mathcal{L}_A$, and $\{\mu_1\le\mu_2\le\cdots\le\mu_{2n}\}$ represent the eigenvalues of $\mathcal{L}_S$. From Lemma \ref{Hua}, the normalized Laplacian spectrum of the graph $P_n$ is given by $$\{0=\rho_1<\rho_2\le\cdots\le\rho_{3n}\}\bigcup \{\mu_1\le\mu_2\le\cdots\le\mu_{2n}\}.$$

Let $M$ be a $p\times q$ matrix and $Y\subset \{1,2,\ldots,p\}$ and $Z\subset \{1,2,\ldots,q\}$. By $M(Y|Z)$, we correspond to the submatrix of $M$ that is formed by eliminating the corresponding rows and columns in $Y$ and $Z$, respectively. We denote the cofactor of $w_{kl}$ of a square matrix $W=(w_{ij})$ of order $p$ by $\operatorname{cf}(W(k,l)).$ We refer to a column vector of appropriate size $e_i,$ with all zeros except the $i$-th element, which is $1$. The (classical) adjoint matrix $\operatorname{adj}(W)$ of a square matrix $W=(w_{ij})$ is defined by   $\operatorname{adj}(W)=(\operatorname{cf}(W(i,j)))^T.$

\section{The derivation of degree-Kirchhoff index}

By Theorem \ref{degKiralt}, the degree-Kirchhoff index for $P_n,$ $(n\geq 2)$  is \begin{eqnarray}
\operatorname{Kf}^*(P_n)&=&14n\left(\sum\limits_{i=2}^{3n} \frac{1}{\rho_i}+\sum\limits_{j=1}^{2n} \frac{1}{\mu_j}\right)\label{kiraltdef1}.
\end{eqnarray}

We obtain from the connection between roots and coefficients $\sum\limits_{i=2}^{3n} \frac{1}{\rho_i}=-\frac{\gamma_{3n-2}}{\gamma_{3n-1}}$ and $\sum\limits_{j=1}^{2n} \frac{1}{\mu_i}=-\frac{\delta_{2n-1}}{\delta_{2n}}=-\frac{\delta_{2n-1}}{\det(\mathcal{L}_S)}\cdot$ 

Hence, from (\ref{kiraltdef1}), for $n\geq 2$, \begin{eqnarray}\operatorname{Kf}^*(P_n)=14n\left(-\frac{\gamma_{3n-2}}{\gamma_{3n-1}}-\frac{\delta_{2n-1}}{\det(\mathcal{L}_S)}\right)\cdot\label{kiraltdef1.1}
\end{eqnarray}

Similarly, for $n\geq 2$, \begin{eqnarray}\operatorname{Kf}^*(P_n')&=&14n\left(\sum\limits_{i=2}^{3n} \frac{1}{\rho_i}+\sum\limits_{j=1}^{2n} \frac{1}{\mu_j'}\right)\notag\\
&=&14n\left(-\frac{\gamma_{3n-2}}{\gamma_{3n-1}}-\frac{\delta'_{2n-1}}{\det(\mathcal{L}_S')}\right)\label{kiraltdef1.2}.
\end{eqnarray}

To evaluate the ratios involved in the above expressions, the next lemmas are used.

\begin{lemma}[Matrix-determinant lemma]\cite{HJ85}\label{Matdet}
Consider a square matrix $W$ of order $n$. Then $\det(W+uv^T)=\det(W)+v^T\operatorname{adj}(W)u,$ where $u$, and $v$ are column vectors of size $n$.
\end{lemma}

\begin{lemma}[Schur determinant formula]\cite{HJ85}
Let $A=\begin{pmatrix}[c|c]
B_1 & B_2\\\hline
B_3 & B_4
\end{pmatrix}$ be a $(p+q)\times (p+q)$ matrix with the blocks $B_1$, $B_2$, $B_3$ and $B_4$ of the respective orders $p\times p,$ $p\times q,$ $q\times p,$ and $q\times q.$ If $B_1$ is an invertible martix, then $\det(A)=\det(B_1)\cdot \det\left(B_4-B_3B_1^{-1}B_2\right).$  
\end{lemma}

In \cite{Sah22}, Sahir and Nayeem have deduced the following lemmas which we will use in our work also.

\begin{lemma}{\cite{Sah22}}\label{det}
Let $R_i=\begin{pmatrix}
-2 &1&0& \cdots&0&0&0\\
1 &-2&1& \cdots&0&0&0\\
0 & 1&-2& \cdots &0&0&0\\
\vdots & \vdots & \vdots&\ddots   & \vdots&\vdots& \vdots\\
0 &0&0& \cdots&-2&1&0\\
0 &0&0& \cdots&1&-2&1\\
0 &0&0& \cdots&0&1&-2
\end{pmatrix}$ be a $i\times i$ matrix. Then $\det(R_i)=(-1)^i\cdot(1+i).$
\end{lemma}

\begin{lemma}{\cite{Sah22}}\label{det1.1}
Let $$R_{i,j}= 
\begin{pmatrix}[cccc|c|cccc]
-2 &1&\cdots&0&0& 0&\cdots&0&0\\
1&-2 &\cdots&0&0& 0&\cdots&0&0\\
0 &1&\cdots&0&0& 0&\cdots&0&0\\
\vdots & \vdots& \ddots &\vdots   & \vdots&  \vdots& \ddots& \vdots & \vdots  \\
0&0&\cdots&-2 &1& 0&\cdots&0&0\\\hline 
0&0&\cdots&1 &-3& 1&\cdots&0&0\\\hline
0&0&\cdots&0 &1& -2&\cdots&0&0\\
\vdots & \vdots& \ddots & \vdots   & \vdots&  \vdots& \ddots& \vdots & \vdots  \\
0&0&\cdots&0 &0& 0&\cdots&-2&1\\
0&0&\cdots&0 &0& 0&\cdots&1&-2\\
\end{pmatrix}$$ be a $i\times i$ matrix where the $(j,j)$-th $(j\leq i)$ element is $-3$. Then $\det(R_{i,j})=(-1)^i\cdot(1+i+j+ij-j^2)$.
\end{lemma}

\begin{lemma}{\cite{Sah22}}\label{lemma1}
 Let $$R=\begin{pmatrix}[ccccc|c|cccc]-2 &1&0&\cdots&0&0& 0&\cdots&0&1\\
1&-2 &1&\cdots&0&0& 0&\cdots&0&0\\
0 &1&-2&\cdots&0&0& 0&\cdots&0&0\\
\vdots & \vdots & \vdots& \ddots&\vdots   & \vdots&  \vdots& \ddots& \vdots & \vdots  \\
0&0&0&\cdots&-2 &1& 0&\cdots&0&0\\\hline 
0&0&0&\cdots&1 &-3& 1&\cdots&0&0\\\hline
0&0&0&\cdots&0 &1& -2&\cdots&0&0\\
\vdots & \vdots & \vdots& \ddots&\vdots   & \vdots&  \vdots& \ddots& \vdots & \vdots  \\
0&0&0&\cdots&0 &0& 0&\cdots&-2&1\\
1&0&0&\cdots&0 &0& 0&\cdots&1&-2\\
\end{pmatrix}$$ be a $2n\times 2n$ matrix where the $(2i,2i)$-th element of $R$ is $-3$. Then $\det(R)=2n.$ 
\end{lemma}

\begin{lemma}{\cite{Sah22}}\label{lemma2}
 Let $$S=\begin{pmatrix}[cccccc|cccccc]
 -2 &1&0&\cdots&0&0&0& 0&\cdots&0&1\\
1&-2 &1&\cdots&0&0&0& 0&\cdots&0&0\\
0 &1&-2&\cdots&0&0&0& 0&\cdots&0&0\\
\vdots&\vdots & \vdots& \ddots & \vdots&\vdots   & \vdots&  \vdots& \ddots& \vdots & \vdots  \\
0&0&0&\cdots&-2&1 &0& 0&\cdots&0&0\\
0&0&0&\cdots&1&-2 &0& 0&\cdots&0&0\\\hline
0&0&0&\cdots&0&0 &-2& 1&\cdots&0&0\\
0&0&0&\cdots&0&0 &1& -2&\cdots&0&0\\
\vdots&\vdots & \vdots& \ddots & \vdots&\vdots   & \vdots&  \vdots& \ddots& \vdots & \vdots  \\
0&0&0&\cdots&0&0 &0& 0&\cdots&-2&1\\
1&0&0&\cdots&0&0 &0& 0&\cdots&1&-2\\
\end{pmatrix}$$ be a $(2n-1)\times (2n-1)$ matrix where the respective size of the diagonal blocks are  $(i-n-1)$ and $(3n-i),$ $i \in \{n+1,n+2,\ldots,3n\}$. Then $\det(S)=-2n.$  
\end{lemma}

\begin{lemma}\label{lemma3}
  Let $$W_1=
\begin{pmatrix}[ccccc|c|ccc|c|cccc]-2 &1&0&\cdots&0&0& 0&\cdots&0&0& 0&\cdots&0&1\\
1&-2 &1&\cdots&0&0& 0&\cdots&0&0& 0&\cdots&0&0\\
0 &1&-2&\cdots&0&0& 0&\cdots&0&0& 0&\cdots&0&0\\
\vdots & \vdots & \vdots& \ddots&\vdots   & \vdots&  \vdots& \ddots&\vdots   & \vdots&  \vdots& \ddots& \vdots & \vdots  \\
0&0&0&\cdots&-2 &1& 0&\cdots&0&0& 0&\cdots&0&0\\\hline 
0&0&0&\cdots&1 &-3& 1&\cdots&0&\cdots& 0&\cdots&0&0\\\hline
0&0&0&\cdots&0 &1& -2&\cdots&0&0& 0&\cdots&0&0\\
\vdots & \vdots & \vdots& \ddots&\vdots   & \vdots&  \vdots& \ddots&\vdots   & \vdots&  \vdots& \ddots& \vdots & \vdots  \\
0&0&0&\cdots&0&0& 0&\cdots&-2 &1& 0&\cdots&0&0\\\hline 
0&0&0&\cdots&0&\cdots& 0&\cdots&1 &-3& 1&\cdots&0&0\\\hline
0&0&0&\cdots&0&0& 0&\cdots&0 &1& -2&\cdots&0&0\\
\vdots & \vdots & \vdots& \ddots&\vdots   & \vdots&  \vdots& \ddots&\vdots   & \vdots&  \vdots& \ddots& \vdots & \vdots  \\
0&0&0&\cdots&0 &0& 0&\cdots&0&0& 0&\cdots&-2&1\\
1&0&0&\cdots&0 &0& 0&\cdots&0&0& 0&\cdots&1&-2\\
\end{pmatrix}$$ be a $2n\times 2n$ matrix with the $(2i, 2i)$-th element and $(2j, 2j)$-th element both $-3$, where $i,j(>i) \in \{1,2,\ldots,n\}$. Then $\det(W_1)=4n+8ij-4n(i-j)-4(i^2+j^2).$
\end{lemma}
\noindent\textit{Proof.} 
Let $$T_{p,s,t}=R_{p,s}-e_te_t^T=\begin{pmatrix}[ccccc|c|ccc|c|cccc]-2 &1&0&\cdots&0&0& 0&\cdots&0&0& 0&\cdots&0&0\\
1&-2 &1&\cdots&0&0& 0&\cdots&0&0& 0&\cdots&0&0\\
0 &1&-2&\cdots&0&0& 0&\cdots&0&0& 0&\cdots&0&0\\
\vdots & \vdots & \vdots& \ddots&\vdots   & \vdots&  \vdots& \ddots&\vdots   & \vdots&  \vdots& \ddots& \vdots & \vdots  \\
0&0&0&\cdots&-2 &1& 0&\cdots&0&0& 0&\cdots&0&0\\\hline 
0&0&0&\cdots&1 &-3& 1&\cdots&0&\cdots& 0&\cdots&0&0\\\hline
0&0&0&\cdots&0 &1& -2&\cdots&0&0& 0&\cdots&0&0\\
\vdots & \vdots & \vdots& \ddots&\vdots   & \vdots&  \vdots& \ddots&\vdots   & \vdots&  \vdots& \ddots& \vdots & \vdots  \\
0&0&0&\cdots&0&0& 0&\cdots&-2 &1& 0&\cdots&0&0\\\hline 
0&0&0&\cdots&0&\cdots& 0&\cdots&1 &-3& 1&\cdots&0&0\\\hline
0&0&0&\cdots&0&0& 0&\cdots&0 &1& -2&\cdots&0&0\\
\vdots & \vdots & \vdots& \ddots&\vdots   & \vdots&  \vdots& \ddots&\vdots   & \vdots&  \vdots& \ddots& \vdots & \vdots  \\
0&0&0&\cdots&0 &0& 0&\cdots&0&0& 0&\cdots&-2&1\\
0&0&0&\cdots&0 &0& 0&\cdots&0&0& 0&\cdots&1&-2\\
\end{pmatrix}$$ be a $p\times p$ matrix where each of the $(s, s)$-th and $(t, t)$-th elements of $T_{p,s,t}$ is $-3,$ $s,t(>s) \in \{1,2,\ldots,p\}$. Apply Lemma \ref{Matdet} to get,
\begin{eqnarray}
\det(T_{p,s,t})&=&\det(R_{p,s})-\operatorname{cf}(R_{p,s}(t,t))\notag\\
&=& \det(R_{p,s})-(-1)^{t+t}\det(R_{t-1,s})\cdot\det(R_{p-t})\notag\\
&=& (-1)^p\cdot(1+p+s+ps-s^2)\notag\\
&&{~~~~~~~~~~~~~~~~~~~~~~~~~~~}-(-1)^{t-1}\cdot(1+(t-1)+s+(t-1)s-s^2)\cdot (-1)^{p-t}(1+p-t)\notag\\
&& \text{~~~~~~~~~~~~~~~~~~~~~~~~~~~~~~~~~~~~~~~~~~~~~~~~~~~~~~~~~~~~~~~~~~~~~~~~~~~~~~~~~~~~~~~~~~~~~(by Lemma \ref{det} and Lemma \ref{det1.1})}\notag\\
&=& (-1)^p\cdot\left(pt(s+1)+s(p+t)+(p+s+t+1)-s^2(p+2-t)-t^2(s+1)\right). \label{3parametersdet}
\end{eqnarray} 
Now, applying  Lemma \ref{Matdet} two consecutive times and using  (\ref{3parametersdet}), we obtain
\begin{eqnarray*}
\det(W_1)&=&\det(T_{2n,2i,2j})-\det(T_{2n-2,2i-1,2j-1})-2\\
&=&4n+8ij-4n(i-j)-4(i^2+j^2).
\end{eqnarray*}

\qed 

\begin{lemma}\label{lemma4}
Let $$W_2=\begin{pmatrix}
R_{i-n-1} & \mathbf{0} & e_{1}e_{3n-j}^T\\
\mathbf{0} & R_{j-i-1}&\mathbf{0}\\
e_{3n-j}e_{1}^T & \mathbf{0}& R_{3n-j}
\end{pmatrix}$$ with $i,j(>i) \in \{n+1, n+2, \ldots, 3n\}$, be a $(2n-2)\times (2n-2)$ matrix. Then $\det(W_2)=2ij+2n(j-i)-(i^2+j^2).$ 
\end{lemma}
\noindent\textit{Proof.} Apply Lemma \ref{Matdet} two times successively to get,
\begin{eqnarray*}
\det(W_2)&=&\det(R_{i-n-1})\cdot\det(R_{j-i-1})\cdot\det(R_{3n-j})-\det(R_{i-n-2})\cdot\det(R_{j-i-1})\cdot(R_{3n-j-1})\\
&=& \det(R_{j-i-1})\cdot\left(\det(R_{i-n-1})\cdot\det(R_{3n-j})-\det(R_{i-n-2})\cdot \det(R_{3n-j-1})\right)\\
&=& (-1)^{j-i-1}(j-i)\\
&&\cdot\left((-1)^{i-n-1}(i-n)\cdot (-1)^{3n-j}(3n-j+1)-(-1)^{i-n-2}(i-n-1)\cdot (-1)^{3n-j-1}(3n-j)\right)\\
&=& (j-i)\cdot\left((i-n)\cdot (3n-j+1)-(i-n-1)\cdot (3n-j)\right)\\
&=& 2ij+2n(j-i)-(i^2+j^2).
\end{eqnarray*}\qed

\begin{lemma}\label{lemma5}
Let $$P=\begin{pmatrix}[ccccc|c|cccc]-2 &1&0&\cdots&0&0& 0&\cdots&0&1\\
	1&-2 &1&\cdots&0&0& 0&\cdots&0&0\\
	0 &1&-2&\cdots&0&0& 0&\cdots&0&0\\
	\vdots & \vdots & \vdots& \ddots&\vdots   & \vdots&  \vdots& \ddots& \vdots & \vdots  \\
	0&0&0&\cdots&-2 &1& 0&\cdots&0&0\\\hline 
	0&0&0&\cdots&1 &-3& 1&\cdots&0&0\\\hline
	0&0&0&\cdots&0 &1& -2&\cdots&0&0\\
	\vdots & \vdots & \vdots& \ddots&\vdots   & \vdots&  \vdots& \ddots& \vdots & \vdots  \\
	0&0&0&\cdots&0 &0& 0&\cdots&-2&1\\
	1&0&0&\cdots&0 &0& 0&\cdots&1&-2\\
\end{pmatrix}$$ be a ${2n\times 2n}$ matrix where the $(2i,2i)$-th element is $-3$, $W_3=P(\{j-n\}|\{j-n\}),$ and $i \in \{1,2,\ldots,n\}$,  $j \in \{n+1,n+2,\ldots,3n\}$. Then
$$\det(W_3)=\begin{cases}
			(j-n)^2-4(j-n)i+4i^2-2n-2n(j-n-2i), & \text{if~ $2i\leq j-n$}\\
            	(j-n)^2-4(j-n)i+4i^2-2n+2n(j-n-2i), & \text{if~ $2i> j-n$ }.
		 \end{cases}$$
\end{lemma}
\noindent\textit{Proof.} Apply Lemma \ref{Matdet} for the cases below.

\begin{case} \mbox{When $2i\leq j-n.$}

$\det(W_3) =\det(R_{j-n-1,2i})\cdot\det(R_{3n-j})-\det(R_{j-n-2,2i-1})\cdot\det(R_{3n-j-1})$\\
$= (-1)^{j-n-1}(1+(j-n-1)+(2i)+(j-n-1)(2i)-(2i)^2)\cdot(-1)^{3n-j}(1+3n-j)$\\
$-(-1)^{j-n-2}(1+(j-n-2)+(2i-1)+(j-n-2)(2i-1)-(2i-1)^2)\cdot(-1)^{3n-j-1}(3n-j)$ (by Lemmas \ref{det} and \ref{det1.1})\\
$=-[(1+(j-n-1)+(2i)+(j-n-1)(2i)-(2i)^2)\cdot(1+3n-j)-(1+(j-n-2)+(2i-1)+(j-n-2)(2i-1)-(2i-1)^2)\cdot(3n-j)]\\
=(j-n)^2-4(j-n)i+4i^2-2n-2n(j-n-2i).$ 
\end{case}

\begin{case} \mbox{When $2i> j-n.$}

$\det(W_3) =\det(R_{j-n-1})\det(R_{3n-j,2n-2i+1})-\det(R_{j-n-2})\det(R_{3n-j-1,2n-2i})$\\
$= (-1)^{j-n-1}(j-n)\cdot(-1)^{3n-j}(1+(3n-j)+(2n-2i+1)+(3n-j)(2n-2i+1)-(2n-2i+1)^2)$\\
$-(-1)^{j-n-2}(j-n-1)\cdot(-1)^{3n-j-1}(1+(3n-j-1)+(2n-2i)+(3n-j-1)(2n-2i)-(2n-2i)^2)$ (by Lemmas \ref{det} and \ref{det1.1})\\
$=-[(j-n)\cdot(1+(3n-j)+(2n-2i+1)+(3n-j)(2n-2i+1)-(2n-2i+1)^2)-(j-n-1)(1+(3n-j-1)+(2n-2i)+(3n-j-1)(2n-2i)-(2n-2i)^2)]\\
=(j-n)^2-4(j-n)i+4i^2-2n+2n(j-n-2i).$ 
\end{case}
\qed 

\begin{lemma}{\cite{Sah22}}\label{coefficientofbeta}
 Let $$N=\begin{pmatrix}
4 &-1&0&0& \cdots&0&0&-1\\
-1&3 &-1&0& \cdots&0&0&0\\
0 &-1&4&-1& \cdots&0&0&0\\
0&0&-1&3 & \cdots&0&0&0\\
\vdots & \vdots & \vdots   & \vdots& \ddots& \vdots & \vdots & \vdots \\
0&0&0&0 & \cdots&3&-1&0\\
0&0&0&0 & \cdots&-1&4&-1\\
-1&0&0&0 & \cdots&0&-1&3\\
\end{pmatrix}$$
and $$N'=\begin{pmatrix}
4 &-1&0&0& \cdots&0&0&1\\
-1&3 &-1&0& \cdots&0&0&0\\
0 &-1&4&-1& \cdots&0&0&0\\
0&0&-1&3 & \cdots&0&0&0\\
\vdots & \vdots & \vdots   & \vdots& \ddots& \vdots & \vdots & \vdots \\
0&0&0&0 & \cdots&3&-1&0\\
0&0&0&0 & \cdots&-1&4&-1\\
1&0&0&0 & \cdots&0&-1&3\\
\end{pmatrix}$$ be two matrices of order $2n.$ Then
\begin{eqnarray*}
\sum\limits_{i=1}^{2n} \det (N(\{i\}|\{i\}))=\sum\limits_{i=1}^{2n} \det  (N'(\{i\}|\{i\}))=\frac{7\sqrt{6}n}{24}\left[\left(\sqrt{3}+\sqrt{2}\right)^{2n}-\left(\sqrt{3}-\sqrt{2}\right)^{2n}\right],\end{eqnarray*}
$\det(N)=(\sqrt{3}-\sqrt{2})^{2n}+(\sqrt{3}+\sqrt{2})^{2n}-2,$ and $\det(N')=(\sqrt{3}-\sqrt{2})^{2n}+(\sqrt{3}+\sqrt{2})^{2n}+2.$
\end{lemma}

\begin{lemma}\label{vieta1} For $n\geq 2,$ 
$$\sum\limits_{i=2}^{3n} \frac{1}{\rho_i}=-\frac{\gamma_{3n-2}}{\gamma_{3n-1}}=\frac{49n^2+42n-19}{42}\cdot$$
\end{lemma}
\noindent\textit{Proof.} From (\ref{char1}), we have $\gamma_{3n-1}=\sum\limits_{i=1}^{3n} \det(-\mathcal{L}_A(\{i\}|\{i\})).$ 

Now, let $i\in \{1,2,\ldots,n\}$. So,\\ $\det(-\mathcal{L}_A(\{i\}|\{i\}))=\begin{vmatrix}
-I_{n-1}& -\sqrt{2}\mathcal{L}_{01}(\{i\}|\{\}) \\ 
-\sqrt{2}\mathcal{L}_{01}(\{i\}|\{\})^T & Q 
\end{vmatrix}$, where $Q= -\mathcal{L}_{11} - \mathcal{L}_{12}$. Using Schur determinant formula, we get\\ $\det(-\mathcal{L}_A(\{i\}|\{i\}))=\det\left(-I_{n-1}\right)\cdot\det\left(Q+2\mathcal{L}_{01}(\{i\}|\{\})^T\cdot \mathcal{L}_{01}(\{i\}|\{\})\right),$ for $i\in \{1,2,\ldots,n\}$.

Now,
$$Q+2\mathcal{L}_{01}(\{i\}|\{\})^T\cdot \mathcal{L}_{01}(\{i\}|\{\})=\frac{1}{3}R.$$

Hence, from Lemma \ref{lemma1}, $\det(-\mathcal{L}_A(\{i\}|\{i\}))=(-1)^{n-1}\cdot\frac{1}{3^{2n}}\cdot 2n,$ for $i\in \{1,2,\ldots,n\}$.

For $i\in \{n+1,n+2,\ldots,3n\}$,\\ $\det(-\mathcal{L}_A(\{i\}|\{i\}))=\begin{vmatrix}
-I_{n}& -\sqrt{2}\mathcal{L}_{01}(\{\}|\{i-n\}) \\ 
-\sqrt{2}\mathcal{L}_{01}(\{\}|\{i-n\})^T & Q[i-n] 
\end{vmatrix}$. 

Using Schur determinant formula, we have, for $i\in \{n+1,n+2,\ldots,3n\}$,\\ $\det(-\mathcal{L}_A(\{i\}|\{i\}))=\det\left(-I_n\right)\cdot\det\left(Q[i-n]+2\mathcal{L}_{01}(\{\}|\{i-n\})^T\cdot \mathcal{L}_{01}(\{\}|\{i-n\})\right)$.

Now, 
\begin{eqnarray*}Q[i-n]+2\mathcal{L}_{01}(\{\}|\{i-n\})^T\cdot \mathcal{L}_{01}(\{\}|\{i-n\})&=&
\frac{1}{3}S.
\end{eqnarray*}

Thus from Lemma \ref{lemma2}, $\det(-\mathcal{L}_A(\{i\}|\{i\}))=(-1)^{n+1}\cdot\frac{1}{3^{2n-1}}\cdot2n,$ where $i\in \{n+1,n+2,\ldots,3n\}$. 

Therefore, 
\begin{eqnarray}
\nonumber     \gamma_{3n-1}&=&\sum\limits_{i=1}^{3n} \det(-\mathcal{L}_A(\{i\}|\{i\}))\\
\nonumber     &=&\sum\limits_{i=1}^{n} \det(-\mathcal{L}_A(\{i\}|\{i\}))+\sum\limits_{i=n+1}^{3n} \det(-\mathcal{L}_A(\{i\}|\{i\}))\\
\nonumber     &=&(-1)^{n-1}\cdot\frac{1}{3^{2n}}\cdot 2n^2+(-1)^{n+1}\cdot\frac{1}{3^{2n-1}}\cdot4n^2\\
&=&(-1)^{n-1}\cdot \frac{14n^2}{3^{2n}}\cdot
\label{a3n-1}
\end{eqnarray}

Again, from (\ref{char1}) 
\begin{eqnarray*}\gamma_{3n-2}&=&\sum\limits_{1\leq i <j\leq 3n} \det(-\mathcal{L}_A(\{i,j\}|\{i,j\}))\\
&=&\sum\limits_{1\leq i <j\leq n} \det(-\mathcal{L}_A(\{i,j\}|\{i,j\}))+\sum\limits_{n+1\leq i <j\leq 3n} \det(-\mathcal{L}_A(\{i,j\}|\{i,j\}))\\
&&+\sum\limits_{\substack{1\leq i \leq n,\\ n+1\leq j\leq 3n}} \det(-\mathcal{L}_A(\{i,j\}|\{i,j\})).\end{eqnarray*}

For $i,j(>i) \in \{1,2,\ldots,n\}$, using Schur determinant formula,  
\begin{eqnarray*}
\det(-\mathcal{L}_A(\{i,j\}|\{i,j\}))&=&\begin{vmatrix}
-I_{n-2}& -\sqrt{2}\mathcal{L}_{01}(\{i,j\}|\{\}) \\ 
-\sqrt{2}\mathcal{L}_{01}(\{i,j\}|\{\})^T & Q 
\end{vmatrix}\\
&=&\det(-I_{n-2})\cdot \det(Q+2\mathcal{L}_{01}(\{i,j\}|\{\})^T\cdot \mathcal{L}_{01}(\{i,j\}|\{\}))\\
&=& (-1)^{n-2}\cdot  \det\left(\frac{1}{3}W_1\right),
\end{eqnarray*}
where $Q+2\mathcal{L}_{01}(\{i,j\}|\{\})^T\cdot \mathcal{L}_{01}(\{i,j\}|\{\})=\frac{1}{3}W_1.$  

By Lemma \ref{lemma3}, 
$\det(W_1)=4n+8ij-4n(i-j)-4(i^2+j^2)$ and hence,\\ $\det(-\mathcal{L}_A(\{i,j\}|\{i,j\}))=(-1)^{n-2}\cdot\frac{1}{3^{2n}}\cdot4\left[n+2ij-n(i-j)-(i^2+j^2)\right],$ for $i,j(>i) \in \{1,2,\ldots,n\}$.

For $i,j(>i) \in \{n+1,n+2,\ldots,3n\},$ using Schur determinant formula,\\  
\begin{eqnarray*}
\det(-\mathcal{L}_A(\{i,j\}|\{i,j\}))&=&\begin{vmatrix}
-I_{n}& -\sqrt{2}\mathcal{L}_{01}(\{\}|\{i-n,j-n\}) \\ 
-\sqrt{2}\mathcal{L}_{01}(\{\}|\{i-n,j-n\})^T & Q(\{i-n,j-n\}|\{i-n,j-n\})
\end{vmatrix}\end{eqnarray*}
\begin{eqnarray*}=&&\det(-I_{n}) \cdot \det(Q(\{i-n,j-n\}|\{i-n,j-n\})+2\mathcal{L}_{01}(\{\}|\{i-n,j-n\})^T\cdot\mathcal{L}_{01}(\{\}|\{i-n,j-n\}))\\
=&&(-1)^n\cdot \det\left(\frac{1}{3}W_2\right),\end{eqnarray*}
where $Q(\{i-n,j-n\}|\{i-n,j-n\})+2\mathcal{L}_{01}(\{\}|\{i-n,j-n\})^T\cdot \mathcal{L}_{01}(\{\}|\{i-n,j-n\})=\frac{1}{3}W_2.$

By Lemma \ref{lemma4},  $\det(W_2)=2ij+2n(j-i)-(i^2+j^2)$, hence, $$\det(-\mathcal{L}_A(\{i,j\}|\{i,j\}))=(-1)^n\cdot\frac{1}{3^{2n-2}}\cdot[2ij+2n(j-i)-(i^2+j^2)], \mbox{ for } i,j(>i) \in \{n+1,n+2,\ldots,3n\}.$$
 
In a similar way, for $i \in \{1,2,\ldots,n\},$ and $j \in \{n+1,n+2,\ldots,3n\},$  
\begin{eqnarray*}
\det(-\mathcal{L}_A(\{i,j\}|\{i,j\}))&=&\begin{vmatrix}
-I_{n-1}& -\sqrt{2}\mathcal{L}_{01}(\{i\}|\{j-n\}) \\ 
-\sqrt{2}\mathcal{L}_{01}(\{i\}|\{j-n\})^T & Q(\{j-n\}|\{j-n\})
\end{vmatrix}\end{eqnarray*}
\begin{eqnarray*}=&&\det(-I_{n}) \cdot \det(Q(\{j-n\}|\{j-n\})+2\mathcal{L}_{01}(\{i\}|\{j-n\})^T\cdot \mathcal{L}_{01}(\{i\}|\{j-n\}))\\
=&&(-1)^{n-1}\cdot \det\left(\frac{1}{3}W_3\right),\end{eqnarray*}
where $Q(\{j-n\}|\{j-n\})+2\mathcal{L}_{01}(\{i\}|\{j-n\})^T\cdot \mathcal{L}_{01}(\{i\}|\{j-n\})=\frac{1}{3}W_3.$
From Lemma \ref{lemma5}, $$\det(W_3)=\begin{cases}
			(j-n)^2-4(j-n)i+4i^2-2n-2n(j-n-2i), & \text{if~ $2i\leq j-n$}\\
            	(j-n)^2-4(j-n)i+4i^2-2n+2n(j-n-2i), & \text{if~ $2i> j-n$ }.
		 \end{cases}$$
Hence,
\begin{eqnarray*}
\gamma_{3n-2}&=&(-1)^{n-2}\cdot\frac{1}{3^{2n+1}}\cdot (n^4+6n^3-7n^2)+(-1)^{n}\cdot\frac{1}{3^{2n-1}}\cdot (4n^4-n^2)+(-1)^n\cdot\frac{1}{3^{2n}}\cdot (4n^4+12n^3-n^2)\\
&=&(-1)^n\cdot \frac{49n^2+42n-19n^2}{3^{2n+1}}\cdot 
\end{eqnarray*}

Hence $\sum\limits_{i=2}^{3n} \frac{1}{\rho_i}=-\frac{\gamma_{3n-2}}{\gamma_{3n-1}}=\frac{49n^2+42n-19}{42}\cdot$ \qed

\begin{lemma}\label{vieta2} For $n\geq 2,$ \\
 $$\sum\limits_{i=1}^{2n} \frac{1}{\mu_i}=\frac{\frac{7\sqrt{6}n}{8}\left[\left(\sqrt{3}+\sqrt{2}\right)^{2n}-\left(\sqrt{3}-\sqrt{2}\right)^{2n}\right]}{(\sqrt{3}+\sqrt{2})^{2n}+(\sqrt{3}-\sqrt{2})^{2n}-2}$$ \\
and \\

$$\sum\limits_{i=1}^{2n} \frac{1}{\mu'_i}=\frac{\frac{7\sqrt{6}n}{8}\left[\left(\sqrt{3}+\sqrt{2}\right)^{2n}-\left(\sqrt{3}-\sqrt{2}\right)^{2n}\right]}{(\sqrt{3}+\sqrt{2})^{2n}+(\sqrt{3}-\sqrt{2})^{2n}+2}.$$
\end{lemma}
\noindent\textit{Proof.} We have $\delta_{2n-1}=\delta'_{2n-1}=(-1)^{2n-1}\cdot\sum\limits_{i=1}^{2n}\mathcal{L}_S(\{i\}|\{i\})=(-1)^{2n-1}\cdot\sum\limits_{i=1}^{2n}\mathcal{L}_S'(\{i\}|\{i\}).$
Since
$\mathcal{L}_S=\frac{1}{3}N$ and $\mathcal{L}'_S=\frac{1}{3}N',$ \\
so, from Lemma \ref{coefficientofbeta}
\begin{eqnarray*}
&&-\delta_{2n-1}=-\delta'_{2n-1}=\frac{1}{3^{2n-1}}\sum\limits_{i=1}^{2n} \det(N(\{i\}|\{i\}))=\frac{1}{3^{2n-1}}\sum\limits_{i=1}^{2n} \det(N'(\{i\}|\{i\}))\\
&=&\frac{1}{3^{2n-1}}\left(\frac{7\sqrt{6}n}{24}\left[\left(\sqrt{3}+\sqrt{2}\right)^{2n}-\left(\sqrt{3}-\sqrt{2}\right)^{2n}\right]\right).
\end{eqnarray*}

Hence, from Lemma \ref{coefficientofbeta}, for $n\geq 2,$\\
 \begin{eqnarray*}
 \sum\limits_{i=1}^{2n} \frac{1}{\mu_i}&=&-\frac{\delta_{2n-1}}{\det(\mathcal{L}_S)}\\ &=&\frac{\frac{1}{3^{2n-1}}\left(\frac{7\sqrt{6}n}{24}\left[\left(\sqrt{3}+\sqrt{2}\right)^{2n}-\left(\sqrt{3}-\sqrt{2}\right)^{2n}\right]\right)}{\frac{1}{3^{2n}}\left((\sqrt{3}+\sqrt{2})^{2n}+(\sqrt{3}-\sqrt{2})^{2n}-2\right)}\\
&=& \frac{\frac{7\sqrt{6}n}{8}\left[\left(\sqrt{3}+\sqrt{2}\right)^{2n}-\left(\sqrt{3}-\sqrt{2}\right)^{2n}\right]}{(\sqrt{3}+\sqrt{2})^{2n}+(\sqrt{3}-\sqrt{2})^{2n}-2},
 \end{eqnarray*}
 
and 

 \begin{eqnarray*}
 \sum\limits_{i=1}^{2n} \frac{1}{\mu'_i}&=&-\frac{\delta_{2n-1}}{\det(\mathcal{L}'_S)}\\ &=&\frac{\frac{1}{3^{2n-1}}\left(\frac{7\sqrt{6}n}{24}\left[\left(\sqrt{3}+\sqrt{2}\right)^{2n}-\left(\sqrt{3}-\sqrt{2}\right)^{2n}\right]\right)}{\frac{1}{3^{2n}}\left((\sqrt{3}+\sqrt{2})^{2n}+(\sqrt{3}-\sqrt{2})^{2n}+2\right)}\\
&=& \frac{\frac{7\sqrt{6}n}{8}\left[\left(\sqrt{3}+\sqrt{2}\right)^{2n}-\left(\sqrt{3}-\sqrt{2}\right)^{2n}\right]}{(\sqrt{3}+\sqrt{2})^{2n}+(\sqrt{3}-\sqrt{2})^{2n}+2}.
 \end{eqnarray*}
\qed

Using relations (\ref{kiraltdef1.1}), and (\ref{kiraltdef1.2}) and Lemma \ref{vieta1}, and Lemma \ref{vieta2}, we obtain the following theorem.

\begin{theorem}\label{Kirchhoff}
For $n \geq 2,$ the degree-Kirchhoff index $\operatorname{Kf}^*(P_n)$ of $P_n$ and $\operatorname{Kf}^*(P_n')$ of $P_n'$ are given by,
\begin{eqnarray*}
\operatorname{Kf}^*(P_n)&=&14n\left(\frac{49n^2+42n-19}{42}+ \frac{\frac{7\sqrt{6}n}{8}\left[\left(\sqrt{3}+\sqrt{2}\right)^{2n}-\left(\sqrt{3}-\sqrt{2}\right)^{2n}\right]}{(\sqrt{3}+\sqrt{2})^{2n}+(\sqrt{3}-\sqrt{2})^{2n}-2}\right),    
\end{eqnarray*}

and 

\begin{eqnarray*}
\operatorname{Kf}^*(P_n')&=&14n\left(\frac{49n^2+42n-19}{42}+ \frac{\frac{7\sqrt{6}n}{8}\left[\left(\sqrt{3}+\sqrt{2}\right)^{2n}-\left(\sqrt{3}-\sqrt{2}\right)^{2n}\right]}{(\sqrt{3}+\sqrt{2})^{2n}+(\sqrt{3}-\sqrt{2})^{2n}+2}\right)\cdot    
\end{eqnarray*}

\end{theorem}

The following theorem is derived from Theorem \ref{kemeny} and Theorem \ref{Kirchhoff}.
\begin{theorem}
For $n \geq 2,$ Kemeny's constant $\operatorname{Kc}(P_n)$ of $P_n$ and $\operatorname{Kc}(P_n')$ of $P_n'$ are given by, 
\begin{eqnarray*}
\operatorname{Kc}(P_n)&=&\left(\frac{49n^2+42n-19}{42}+ \frac{\frac{7\sqrt{6}n}{8}\left[\left(\sqrt{3}+\sqrt{2}\right)^{2n}-\left(\sqrt{3}-\sqrt{2}\right)^{2n}\right]}{(\sqrt{3}+\sqrt{2})^{2n}+(\sqrt{3}-\sqrt{2})^{2n}-2}\right),    
\end{eqnarray*}

and 

\begin{eqnarray*}
\operatorname{Kc}(P_n')&=&\left(\frac{49n^2+42n-19}{42}+ \frac{\frac{7\sqrt{6}n}{8}\left[\left(\sqrt{3}+\sqrt{2}\right)^{2n}-\left(\sqrt{3}-\sqrt{2}\right)^{2n}\right]}{(\sqrt{3}+\sqrt{2})^{2n}+(\sqrt{3}-\sqrt{2})^{2n}+2}\right)\cdot    
\end{eqnarray*}

\end{theorem}
We obtain the formulae for the total count of spanning trees of $P_n$ and $P'_n$ using the Theorem \ref{spanningtree} as follows.
\begin{corollary}
The total count of spanning trees for $P_n$ is 
\begin{eqnarray*}\tau(P_n)&=&2^n\cdot n\cdot\left((\sqrt{3}-\sqrt{2})^{2n}+(\sqrt{3}+\sqrt{2})^{2n}-2\right)
\end{eqnarray*}

and the total count of spanning trees for $P'_n$ is \begin{eqnarray*}\tau(P'_n)&=&2^n\cdot n\cdot\left((\sqrt{3}-\sqrt{2})^{2n}+(\sqrt{3}+\sqrt{2})^{2n}+2\right).
\end{eqnarray*}
\end{corollary}

\noindent\textit{Proof.} From Theorem \ref{spanningtree},
\begin{eqnarray*}\tau(P_n)&=&\frac{1}{14n}\cdot\left(\prod\limits_{i=1}^{5n} d_i\right)\cdot\left(\prod\limits_{i=2}^{3n} \rho_i\prod\limits_{j=1}^{2n}\mu_j\right)\\
&=&\frac{2^n3^{4n}\cdot(-1)^{3n-1}\gamma_{3n-1}\cdot \det{(\mathcal{L}_S)}}{14n}\\
&=&\frac{2^n3^{4n}\cdot\frac{14n^2}{3^{2n}}\cdot\frac{1}{3^{2n}}\left((\sqrt{3}-\sqrt{2})^{2n}+(\sqrt{3}+\sqrt{2})^{2n}-2\right)}{14n}\\
	&=&2^n\cdot n\cdot\left((\sqrt{3}-\sqrt{2})^{2n}+(\sqrt{3}+\sqrt{2})^{2n}-2\right). \end{eqnarray*}
and
\begin{eqnarray*}\tau(P_n')&=&\frac{1}{14n}\cdot\left(\prod\limits_{i=1}^{5n} d_i\right)\cdot\left(\prod\limits_{i=2}^{3n} \rho_i\prod\limits_{j=1}^{2n}\mu'_j\right)\\
&=&\frac{2^n3^{4n}\cdot(-1)^{3n-1}\gamma_{3n-1}\cdot \det{(\mathcal{L}_S')}}{14n}\\
&=&\frac{2^n3^{4n}\cdot\frac{14n^2}{3^{2n}}\cdot\frac{1}{3^{2n}}\left((\sqrt{3}-\sqrt{2})^{2n}+(\sqrt{3}+\sqrt{2})^{2n}+2\right)}{14n}\\
	&=&2^n\cdot n\cdot\left((\sqrt{3}-\sqrt{2})^{2n}+(\sqrt{3}+\sqrt{2})^{2n}+2\right). \end{eqnarray*}
	\qed

The same formulae for the total count of spanning trees of $P_n$ and $P_n'$ were obtained by Sahir and Nayeem \cite{Sah22} using the Laplacian eigenvalues of $P_n$ and $P_n'$.

\section{The relations among different indices}
For $P_n$ and $P'_n$, we compute the Schultz and Gutman indices.

\begin{theorem}\label{Gut1}
The Gutman index $\operatorname{Gut}(P_n)$ of the pentagonal cylinder chain $P_n$ on $5n$ vertices is
\[ \operatorname{Gut}(P_n)=\begin{cases} 
      49n^3+64n^2+5n, \textit{~~if~} n\textit{~is~even}\\
      49n^3+64n^2+4n, \textit{~~if~} n\textit{~is~odd.} 
   \end{cases}
\]
\end{theorem}
\noindent\textit{Proof.} 
The vertices of the pentagonal cylinder chain $P_n$ (see Figure \ref{fig1}) can be categorized below.
\begin{enumerate}
  \item An $a$-type vertex is a degree-$3$ vertex that is non-adjacent to any degree-$2$ vertex. So, the vertices with labels $\{1,1',3,3',\ldots,2n-1,2n-1'\}$ are referred to as $a$-type vertices.
  \item A $b$-type vertex is a degree-$3$ vertex that is adjacent to a degree-$2$ vertex. So, the vertices with labels $\{2,2',4,4',\ldots,2n,2n'\}$ are referred to as $b$-type vertices.
  \item A $c$-type vertex is a vertex of degree-$2$. So, the vertices with labels $\{\tilde{1},\tilde{2},\ldots,\tilde{n}\}$ are referred to as $c$-type vertices.
\end{enumerate}
Furthermore, the vertices with labels $\{1,2,3,\ldots,2n\}$ and $\{1',2',3',\ldots,2n'\}$ are said to be upper vertices and lower vertices respectively.  

Now we have the following observations. 

\begin{enumerate}
   \item The aggregate of the weighted distances $d_id_jd_{ij}$ between an $a$-type upper vertex and
  
   \begin{enumerate}
     \item the vertices with labels $\{1,2,3,\ldots,2n\}$ is  $9\left(2\sum\limits_{i=1}^{n-1} i+n\right)=9n^2,$ and
     \item the vertices with labels $\{1',2',3',\ldots,2n'\}$ is $9\left(2\sum\limits_{i=1}^{n} i+n\right)=9(n^2+2n).$
       \end{enumerate}
   
    \item The aggregate of the weighted distances $d_id_jd_{ij}$ between an $a$-type lower vertex and
  
   \begin{enumerate}
     \item the vertices with labels $\{1,2,3,\ldots,2n\}$ is  $9\left(2\sum\limits_{i=1}^{n} i+n\right)=9(n^2+2n),$ and
     \item  the vertices with labels $\{1',2',3',\ldots,2n'\}$ is $9\left(2\sum\limits_{i=1}^{n-1} i+n\right)=9n^2.$
   \end{enumerate}
   
  \item The aggregate of the weighted distances $d_id_jd_{ij}$ between a vertex of type $a$ (regardless of upper or lower) and the vertices with labels $\{\tilde{1},\tilde{2},\ldots,\tilde{n}\}$ is
    $$\begin{cases}
        6\left(4\sum\limits_{i=1}^{n/2} i\right)=6\left\{\frac{n}{2}(n+2)\right\},\mbox{if~} n \mbox{~is even}\\
       6\left(4\sum\limits_{i=1}^{(n+1)/2} i-(n+1)\right)=6\left\{\frac{1}{2}(n+1)^2\right\}, \mbox{if~} n \mbox{~is odd}.
     \end{cases}$$
 
   \item The aggregate of the weighted distances $d_id_jd_{ij}$ between an upper vertex of type $b$ and
  
   \begin{enumerate}
     \item  the vertices with labels $\{1,2,3,\ldots,2n\}$ is  $9\left(2\sum\limits_{i=1}^{n-1} i+n\right)=9n^2,$ and
     \item  the vertices with labels $\{1',2',3',\ldots,2n'\}$ is  $9\left(2\sum\limits_{i=1}^{n} i+(n+1)\right)=9(n+1)^2.$
   \end{enumerate}
   
  \item The aggregate of the weighted distances $d_id_jd_{ij}$ between a lower vertex of type $b$ and
 
   \begin{enumerate}
     \item the vertices with labels $\{1,2,3,\ldots,2n\}$ is  $9\left(2\sum\limits_{i=1}^{n} i+(n+1)\right)=9(n+1)^2,$ and
     \item the vertices with labels $\{1',2',3',\ldots,2n'\}$ is  $9\left(2\sum\limits_{i=1}^{n-1} i+n\right)=9n^2.$
  \end{enumerate}
   
     \item The aggregate of the weighted distances $d_id_jd_{ij}$ between a vertex of type $b$ (regardless of upper or lower) and the vertices with labels $\{\tilde{1},\tilde{2},\ldots,\tilde{n}\}$ is

 $$\begin{cases}
        6\left(2\sum\limits_{i=1}^{n/2} (2i-1)+n\right)=6\left\{\frac{n}{2}(n+2)\right\},\mbox{if~} n \mbox{~is even}\\
       6\left(2\sum\limits_{i=1}^{(n+1)/2} (2i-1)-1\right)=6\left\{\frac{1}{2}(n+1)^2-1\right\}, \mbox{if~} n \mbox{~is odd}.
     \end{cases}$$

 \item The aggregate of the weighted distances $d_id_jd_{ij}$ between a vertex of type $c$ and
  
   \begin{enumerate}
     \item the vertices with labels $\{1,2,3,\ldots,2n\}$ is  $6\left(2\sum\limits_{i=1}^{n} i+n\right)=6(n^2+2n),$
     \item the vertices with labels $\{1',2',3',\ldots,2n'\}$ is  $6\left(2\sum\limits_{i=1}^{n} i+n\right)=6(n^2+2n),$
     \item the vertices with labels $\{\tilde{1},\tilde{2},\ldots,\tilde{n}\}$ is
     $$\begin{cases}
        4\left(4\sum\limits_{i=1}^{n/2} i+n-2\right)=4\left\{\frac{n^2+4n-4}{2}\right\},\mbox{if~} n \mbox{~is even}\\
       4\left(4\sum\limits_{i=1}^{(n+1)/2} i-4\right)=4\left\{\frac{(n+1)(n+3)}{2}-4\right\}, \mbox{if~} n \mbox{~is odd}.
     \end{cases}$$

   \end{enumerate}

 \end{enumerate}

Hence, the Gutman index $\operatorname{Gut}(P_n)$ of the pentagonal cylinder chain $P_n$ is 
\begin{eqnarray*}
\operatorname{Gut}(P_n)&=&\frac{1}{2}\sum\limits_{i=1}^{5n}\sum\limits_{j=1}^{5n} d_id_jd_{ij}\\
&=&\frac{\sum\limits_{j=1}^{2n}\left(9(n^2+n^2+2n+n^2+(n+1)^2)\atop+6(\frac{n}{2}(n+2)+\frac{n}{2}(n+2)+n^2+2n)   \right)+\sum\limits_{j=1}^{n}4\left(\frac{n^2+4n-4}{2}\right)    }{2}
\\&=&\frac{2n\left(48n^2+60n+9\right)+n(2n^2+8n-8)    }{2}\\
&=&49n^3+64n^2+5n, \mbox{if $n$ is even,}
\end{eqnarray*} 
and

\begin{eqnarray*}
\operatorname{Gut}(P_n)&=&\frac{1}{2}\sum\limits_{i=1}^{5n}\sum\limits_{j=1}^{5n} d_id_jd_{ij}\\
&=&\frac{\sum\limits_{j=1}^{2n}\left(9(n^2+n^2+2n+n^2+(n+1)^2)\atop+6(\frac{1}{2}(n+1)^2+\frac{1}{2}(n+1)^2-1+n^2+2n)   \right)+\sum\limits_{j=1}^{n}4\left(\frac{(n+1)(n+3)}{2}-4\right)    }{2}
\\&=&\frac{2n\left(48n^2+60n+9\right)+n(2n^2+8n-10)    }{2}\\
&=&49n^3+64n^2+4n,\mbox{if $n$ is odd.}\end{eqnarray*}  \qed

Next, we present another automorphism of the pentagonal M\"{o}bius chain $P'_n$ below, which will be useful for finding the Gutman index of $P'_n$.

\begin{figure}[h]
\begin{center}
\begin{tikzpicture}[place/.style={thick, circle,draw=black!50,fill=black!20,inner sep=0pt, minimum size = 6 mm}, transform shape, scale=0.8]
\node[place, fill=red!15!] (1) {\footnotesize 1} ;
\node (2)[place, fill=red!15! ] [below of = 1, node distance=2cm] {\footnotesize $1'$} edge[thick] (1);
\node (3)[place, fill=green!15!] [right of = 1, node distance=2cm] {\footnotesize 2} edge[thick] (1);
\node (4)[place, fill=yellow!15!] [below of = 3, node distance=1cm] {\footnotesize $\tilde{1}$} edge[thick] (3) ;
\node (5)[place, fill=green!15!] [ below of = 4] {\footnotesize $2'$} edge[thick] (4) edge[thick] (2);
\node (6)[place] [right of = 3,node distance=2cm] {\footnotesize \tiny $2i-3$} ;
\node (7)[place] [below of = 6,node distance=2cm] {\footnotesize \tiny $2i-3'$} edge[thick] (6) ;
\node (8)[place] [right of = 6,node distance=2cm] {\footnotesize \tiny $2i-2$} edge[thick] (6);
\node (9)[place] [below of = 8, node distance=1cm] {\footnotesize \tiny $\widetilde{i-1}$} edge[thick] (8) ;
\node (10)[place] [right of = 7,node distance=2cm] {\footnotesize \tiny $2i-2'$} edge[thick] (7) edge[thick] (9);
\node (11)[place, fill=red!15!] [right of = 8,node distance=2cm] {\footnotesize \tiny $2i-1$} edge[thick, out=150, in=-30,out looseness=0.8, in looseness=.8] (10);
\node (12)[place, fill=red!15!] [below of = 11,node distance=2cm] {\footnotesize \tiny $2i-1'$} edge[thick, out=-150, in=30,out looseness=0.8, in looseness=.8] (8) edge[thick] (11);
\node (13)[place, fill=green!15!] [right of = 11, node distance=2cm] {\footnotesize  $2i$} edge[thick] (11);
\node (14)[place, fill=yellow!15!] [below of = 13] {\footnotesize $\tilde{i}$} edge[thick] (13);
\node (15)[place, fill=green!15!] [below of = 14, node distance=1cm] {\footnotesize $2i'$}  edge[thick] (14) edge[thick] (12);
\node (16)[place] [right of = 13, node distance=2cm] {\footnotesize \tiny $2n-1$} ;
\node (17)[place] [below of = 16, node distance=2cm] {\footnotesize \tiny $2n-1'$} edge[thick] (16);
\node (18)[place] [right of = 16, node distance=2cm] {\footnotesize $2n$} edge[thick] (16)  edge[thick, out=55, in=125,out looseness=0.4, in looseness=.4] (1);
\node (19)[place] [below of = 18, node distance=1cm] {\footnotesize $\tilde{n}$} edge[thick] (18);
\node (20)[place] [ below of = 19] {\footnotesize $2n'$} edge[thick] (19) edge[thick] (17) edge[thick] (17) edge[thick, out=-55, in=-125, out looseness=0.4, in looseness=0.4] (2);
\draw (3) node [xshift = 1.1cm, yshift =-1 cm]{$\cdots\cdots$};
\draw (13) node [xshift = 1.1cm, yshift =-1 cm]{$\cdots\cdots$};
\end{tikzpicture}
\caption{\label{fig3} Another representation of M\"{o}bius pentagonal chain $P'_n.$}
\end{center}
\end{figure}
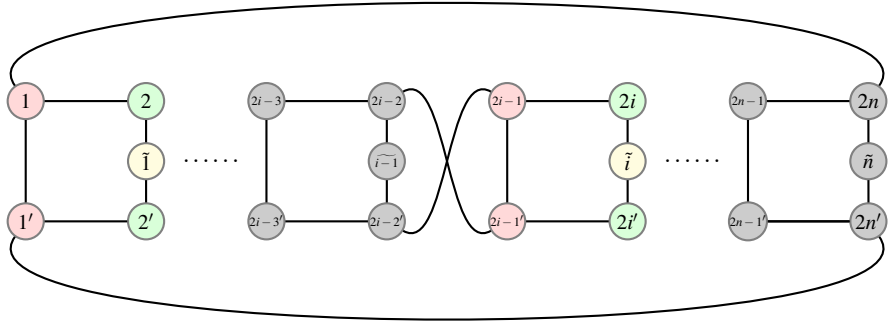

\begin{theorem}\label{Gut2}
The Gutman index $Gut(P'_n)$ of the pentagonal M\"{o}bius chain $P'_n$ on $5n$ vertices is
\[ \operatorname{Gut}(P'_n)=\begin{cases} 
      49n^3+64n^2-13n, \textit{~~if~} n\textit{~is~even}\\
      49n^3+64n^2-14n, \textit{~~if~} n\textit{~is~odd.} 
   \end{cases}
\]
\end{theorem}
\noindent\textit{Proof.} 
 We call the subgraph induced by the vertices with labels $\{(2i-1),(2i-1)',2i,2i',\tilde{i}\}$ as the $i$-th block of $P_n'$.

Let $\pi^k:V(P_n') \rightarrow V(P_n')$ be defined by $\pi^k(i)\equiv (i+2k) \pmod{2n}$, $\pi^k(i')\equiv \left((i+2k) \pmod{2n}\right)',$ and $\pi^k(\tilde{i})=\tilde{r}$ where  $r\equiv (i+k) \pmod{n}$ for each $k\in\{1,2,\ldots, n-1\}$. Then for each $k\in\{1,2,\ldots, n-1\},\pi^k$ is an automorphism of the pentagonal M\"{o}bius chain $P_n'$ (see Figure \ref{fig3}), and hence the aggregate of the weighted distances $d_id_jd_{ij}$ from the vertices in each block to all the vertices of $P_n'$ is same. So, we calculate the aggregate of the weighted distances $d_id_jd_{ij}$ from the vertices in the first block to all the vertices of $P_n'$, and then multiply that quantity by $n$ to get the sum of the weighted distances $d_id_jd_{ij}$ between any pair of vertices of $P_n'.$

\begin{enumerate}
	\item The aggregate of the weighted distances $d_id_jd_{ij}$ between vertex with label $1$ and
  
   \begin{enumerate}
     \item the vertices with labels $\{1,2,3,\ldots,2n\}$ is  $9\left(2\sum\limits_{i=1}^{n} i-1\right)=9\left\{n(n+1)-1\right\},$
	\item the vertices with labels $\{1',2',3',\ldots,2n'\}$ is $9\left(2\sum\limits_{i=1}^{n} i\right)=9\left\{n(n+1)\right\}.$
	\end{enumerate}
	
	\item The aggregate of the weighted distances $d_id_jd_{ij}$ between vertex with label $1'$ and
  
   \begin{enumerate}
     \item the vertices with labels $\{1,2,3,\ldots,2n\}$ is  $9\left(2\sum\limits_{i=1}^{n} i\right)=9\left\{n(n+1)\right\},$
	\item the vertices with labels $\{1',2',3',\ldots,2n'\}$ is $9\left(2\sum\limits_{i=1}^{n} i-1\right)=9\left\{n(n+1)-1\right\}.$
	\end{enumerate}
    \item The aggregate of the weighted distances $d_id_jd_{ij}$ between vertex with label $1$ or with label $1'$ and the vertices with labels $\{\tilde{1},\tilde{2},\ldots,\tilde{n}\}$ is
	$$\begin{cases}
       6\left(4\sum\limits_{i=1}^{n/2} i\right)=6\left\{\frac{n}{2}(n+2)\right\},\mbox{if~} n \mbox{~is even}\\
       6\left(4\sum\limits_{i=1}^{(n+1)/2} i-(n+1)\right)=6\left\{\frac{1}{2}(n+1)^2\right\}, \mbox{if~} n \mbox{~is odd}.
     \end{cases}$$
	
\item The aggregate of the weighted distances $d_id_jd_{ij}$ between vertex with label $2$ and
  
   \begin{enumerate}
     \item the vertices with labels $\{1,2,3,\ldots,2n\}$ is  $9\left(2\sum\limits_{i=1}^{n} i-2\right)=9\left\{n(n+1)-2\right\},$
	 \item the vertices with labels $\{1',2',3',\ldots,2n'\}$ is  $9\left(2\sum\limits_{i=1}^{n} i+2\right)=9\left\{n(n+1)+2\right\}.$
	\end{enumerate}
	
	\item The aggregate of the weighted distances $d_id_jd_{ij}$ between vertex with label $2'$ and
 
   \begin{enumerate}
     \item the vertices with labels $\{1,2,3,\ldots,2n\}$ is  $9\left(2\sum\limits_{i=1}^{n} i+2\right)=9\left\{n(n+1)+2\right\},$
	 \item the vertices with labels $\{1',2',3',\ldots,2n'\}$ is  $9\left(2\sum\limits_{i=1}^{n} i-2\right)=9\left\{n(n+1)-2\right\}.$
	\end{enumerate}
	
	\item The aggregate of the weighted distances $d_id_jd_{ij}$ between vertex with label $2$ or with label $2'$ and the vertices with labels $\{\tilde{1},\tilde{2},\ldots,\tilde{n}\}$ is
	
	$$\begin{cases}
      6\left(2\sum\limits_{i=1}^{n/2} (2i-1)+n\right)=6\left\{\frac{n}{2}(n+2)\right\},\mbox{if~} n \mbox{~is even}\\
      6\left(2\sum\limits_{i=1}^{(n+1)/2} (2i-1)-1\right)=6\left\{\frac{1}{2}(n+1)^2-1\right\}, \mbox{if~} n \mbox{~is odd}.
     \end{cases}$$

	\item The aggregate of the weighted distances $d_id_jd_{ij}$ between vertex with label $\tilde{1}$ and

   \begin{enumerate}
     \item the vertices with labels $\{1,2,3,\ldots,2n\}$ is  $6\left(2\sum\limits_{i=1}^{n} i+n\right)=6\left\{n^2+2n\right\},$
	 \item the vertices with labels $\{1',2',3',\ldots,2n'\}$ is  $6\left(2\sum\limits_{i=1}^{n} i+n\right)=6\left\{n^2+2n\right\},$
		\item the vertices with labels $\{\tilde{1},\tilde{2},\ldots,\tilde{n}\}$ is
     $$\begin{cases}
        4\left(4\sum\limits_{i=1}^{n/2} i+n-2\right)=4\left\{\frac{n^2+4n-4}{2}\right\},\mbox{if~} n \mbox{~is even}\\
       4\left(4\sum\limits_{i=1}^{(n+1)/2} i-4\right)=4\left\{\frac{(n+1)(n+3)}{2}-4\right\}, \mbox{if~} n \mbox{~is odd}.
     \end{cases}$$
		
	\end{enumerate}

\end{enumerate}

Hence, the Gutman index $\operatorname{Gut}(P'_n)$ of the pentagonal M\"{o}bius chain $P'_n$ is 
\begin{eqnarray*}
\operatorname{Gut}(P'_n)&=&\frac{1}{2}\sum\limits_{i=1}^{5n}\sum\limits_{j=1}^{5n} d_id_jd_{ij}\\
&=&\frac{\sum\limits_{j=1}^{n}\left[2\left(9(n(n+1)+n(n+1)-1+n(n+1)-2+n(n+1)+2)\atop+6(\frac{n}{2}(n+2)+\frac{n}{2}(n+2)+n^2+2n)   \right)+4\left(\frac{n^2+4n-4}{2}\right)  \right]  }{2}
\\&=&\frac{2n\left(48n^2+60n-9\right)+n(2n^2+8n-8)    }{2}\\
&=&49n^3+64n^2-13n, \mbox{~if $n$ is even,}
\end{eqnarray*} 
and

\begin{eqnarray*}
\operatorname{Gut}(P'_n)&=&\frac{1}{2}\sum\limits_{i=1}^{5n}\sum\limits_{j=1}^{5n} d_id_jd_{ij}\\
&=&\frac{\sum\limits_{j=1}^{n}\left[2\left(9(n(n+1)+n(n+1)-1+n(n+1)-2+n(n+1)+2)\atop+6(\frac{1}{2}(n+1)^2+\frac{1}{2}(n+1)^2-1+n^2+2n)   \right)+4\left(\frac{(n+1)(n+3)}{2}-4\right)    \right]}{2}
\\&=&\frac{2n\left(48n^2+60n-9\right)+n(2n^2+8n-10)    }{2}\\
&=&49n^3+64n^2-14n,  \mbox{~if $n$ is odd.}\end{eqnarray*} \qed

The Schultz index for $P_n$ and $P'_n$ are obtained by using a similar procedure.

\begin{theorem}\label{Sc1}
The Schultz index $\operatorname{Sc}(P_n)$ for $P_n$ on $5n$ vertices is
\[ \operatorname{Sc}(P_n)=\begin{cases} 
      35n^3+48n^2+2n, \textit{~~if~} n\textit{~is~even}\\
      35n^3+48n^2+n, \textit{~~if~} n\textit{~is~odd.} 
   \end{cases}
\]
\end{theorem}
\noindent\textit{Proof.} Considering the values of $(d_i+d_j)d_{ij}$ for different cases as described in the proof of Theorem \ref{Gut1}, the Schultz index $\operatorname{Sc}(P_n)$ of the pentagonal cylinder chain $P_n$ is given by
\begin{eqnarray*}
\operatorname{Sc}(P_n)&=&\frac{1}{2}\sum\limits_{i=1}^{5n}\sum\limits_{j=1}^{5n} (d_i+d_j)d_{ij}\\
&=&\frac{\sum\limits_{j=1}^{2n}\left(6(n^2+n^2+2n+n^2+(n+1)^2)\atop+5(\frac{n}{2}(n+2)+\frac{n}{2}(n+2)+n^2+2n)   \right)+\sum\limits_{j=1}^{n}4\left(\frac{n^2+4n-4}{2}\right)    }{2}
\\&=&\frac{2n\left(34n^2+44n+6\right)+n(2n^2+8n-8)    }{2}\\
&=&35n^3+48n^2+2n,   \mbox{~if $n$ is even,}
\end{eqnarray*}
and 
\begin{eqnarray*}
\operatorname{Sc}(P_n)&=&\frac{1}{2}\sum\limits_{i=1}^{5n}\sum\limits_{j=1}^{5n} (d_i+d_j)d_{ij}\\
&=&\frac{\sum\limits_{j=1}^{2n}\left(6(n^2+n^2+2n+n^2+(n+1)^2)\atop+5(\frac{1}{2}(n+1)^2+\frac{1}{2}(n+1)^2-1+n^2+2n)   \right)+\sum\limits_{j=1}^{n}4\left(\frac{(n+1)(n+3)}{2}-4\right)    }{2}
\\&=&\frac{2n\left(34n^2+44n+6\right)+n(2n^2+8n-10)    }{2}\\
&=&35n^3+48n^2+n,  \mbox{~if $n$ is odd.}\end{eqnarray*} \qed

\begin{theorem}\label{Sc2}
The Schultz index $\operatorname{Sc}(P'_n)$ for $P'_n$ on $5n$ vertices is
\[ \operatorname{Sc}(P'_n)=\begin{cases} 
      35n^3+48n^2-10n, \textit{~~if~} n\textit{~is~even}\\
      35n^3+48n^2-11n, \textit{~~if~} n\textit{~is~odd.} 
   \end{cases}
\]
\end{theorem}
\noindent\textit{Proof.}
Similar to Theorem \ref{Gut2}, the Schultz index $\operatorname{Sc}(P'_n)$ of the pentagonal M\"{o}bius chain $P'_n$ is 
\begin{eqnarray*}
\operatorname{Sc}(P'_n)&=&\frac{1}{2}\sum\limits_{i=1}^{5n}\sum\limits_{j=1}^{5n} (d_i+d_j)d_{ij}\\
&=&\frac{\sum\limits_{j=1}^{n}\left[2\left(6(n(n+1)+n(n+1)-1+n(n+1)-2+n(n+1)+2)\atop+5(\frac{n}{2}(n+2)+\frac{n}{2}(n+2)+n^2+2n)   \right)+4\left(\frac{n^2+4n-4}{2}\right) \right]   }{2}
\\&=&\frac{2n\left(34n^2+44n-6\right)+n(2n^2+8n-8)    }{2}\\
&=&35n^3+48n^2-10n,   \mbox{if $n$ is even,}
\end{eqnarray*} 
and

\begin{eqnarray*}
\operatorname{Sc}(P'_n)&=&\frac{1}{2}\sum\limits_{i=1}^{5n}\sum\limits_{j=1}^{5n} (d_i+d_j)d_{ij}\\
&=&\frac{\sum\limits_{j=1}^{n}\left[2\left(6(n(n+1)+n(n+1)-1+n(n+1)-2+n(n+1)+2)\atop+5(\frac{1}{2}(n+1)^2+\frac{1}{2}(n+1)^2-1+n^2+2n)   \right)+4\left(\frac{(n+1)(n+3)}{2}-4\right)   \right] }{2}
\\&=&\frac{2n\left(34n^2+44n-6\right)+n(2n^2+8n-10)    }{2}\\
&=&35n^3+48n^2-11n,   \mbox{if $n$ is odd.}\end{eqnarray*} \qed

Combining the results obtained in Theorem \ref{Gut1}, Theorem \ref{Gut2}, Theorem \ref{Sc1}, and Theorem \ref{Sc2} we have the following relations.
\begin{theorem}
$\operatorname{Gut}(P_n)$ and $\operatorname{Sc}(P_n)$ are related by
$$\operatorname{Gut}(P_n)=(14n^3+16n^2+3n)+\operatorname{Sc}(P_n),$$ and $\operatorname{Gut}(P'_n)$ and $\operatorname{Sc}(P'_n)$ are related by
$$\operatorname{Gut}(P'_n)=(14n^3+16n^2-3n)+\operatorname{Sc}(P'_n).$$
\end{theorem}

There are no explicit relations between the degree-Kirchhoff index $\operatorname{Kf}^*(P_n)$ ($\operatorname{Kf}^*(P'_n)$) and the Gutman index $\operatorname{Gut}(P_n)$ ($\operatorname{Gut}(P'_n)$) of the pentagonal cylinder (M\"{obius}) chain $P_n$ ($P'_n$) but when $n$ gets larger, the Gutman index approaches to thrice of the degree-Kirchhoff index. Combining Theorem $\ref{Kirchhoff}$, Theorem $\ref{Gut1}$, and Theorem $\ref{Gut2}$ we get the theorem below. 
\begin{theorem}
$$\lim_{n \to\infty}\frac{\operatorname{Gut}(P_n)}{\operatorname{Kf}^*(P_n)}=3$$
and $$\lim_{n \to\infty}\frac{\operatorname{Gut}(P'_n)}{\operatorname{Kf}^*(P'_n)}=3.$$
\end{theorem}

\begin{table}[h]
\begin{center}
\begin{tabular}{|c|c|c|c|c|c|c|}
\hline
 $n$ &$\operatorname{Gut}(P_n)$&  $\operatorname{Kf}^*(P_n)$ & $\frac{\operatorname{Gut}(P_n)}{\operatorname{Kf}^*(P_n)}$   &$\operatorname{Gut}(P'_n)$&  $\operatorname{Kf}^*(P'_n)$ & $\frac{\operatorname{Gut}(P'_n)}{\operatorname{Kf}^*(P'_n)}$\\ 
 \hline
 2&             658 &296.5&2.2192      & 622 & 291.6 &2.1330\\
 \hline
 3  & 1911 &818.6136 &2.3344    & 1857  & 817.5 & 2.2716\\
 \hline
4 & 4180 &1724.2 &2.4243     &  4108 &1724 & 2.3828\\
 \hline
5 & 7745 &3110.1720 & 2.4902 &  7655 &3110.1404&2.4613\\
 \hline
6 & 12918 & 5074.2273  &2.5458     &  12810 &5074.2227 & 2.5245 \\
\hline
7 & 19971 &7714.3065 & 2.5888  &  19845 & 7714.3059 & 2.5725 \\
 \hline
8 & 29224 &11128.4 & 2.6261        &  29080 &11128.3999 & 2.6131\\
 \hline
9 & 40941 &15414.5062 & 2.6560 &  40779 &15414.5062 &2.6455 \\
 \hline
10 & 55450 &20670.6249 & 2.6826       &  55270 &20670.6249 & 2.6738 \\
 \hline
20 & 417700 &148142.4997 & 2.8196       &  417340 & 148142.4997 & 2.8172\\
\hline
50 & 6285250 & 2151365.6234 & 2.9215   &  6284350 & 2151365.6234 & 2.9211\\
\hline

99 & 48172311 &16278895.2499 & 2.9592  &  48170529 & 16278895.2499 & 2.9591\\
\hline
 \end{tabular}
\caption{\label{tab2} $\operatorname{Gut}(P_n)$ and $\operatorname{Kf}^*(P_n)$ for different $n$.}
\end{center}
\end{table}

\section{Conclusion}
The degree-Kirchhoff index, the Gutman index, and the Schultz index of the pentagonal M\"{o}bius chain $P'_n$ and the pentagonal cylinder chain $P_n$ are computed in closed-form in this study. For the previously mentioned types of graphs, we have deduced a relationship between the Schultz index and the Gutman index. Additionally, it has been discovered that for $P_n$ and $P'_n$, the Gutman index is approximately three times larger than the degree-Kirchhoff index when $n$ is large enough.

\bibliographystyle{abbrv}
\bibliography{mybib.bib}
\end{document}